\documentclass[11pt,reqno]{amsart}

\usepackage{epsfig,amssymb,amsmath,amscd,amsfonts,color, amsthm}
\title[Large deviations for weighted empirical mean]{Large Deviations for weighted empirical mean
with outliers}

\author{M. Ma\"{\i}da, J. Najim
and S. P\'ech\'e}

\date{\today}

\newtheorem{theo}{Theorem}[section]
\newtheorem{lemma}[theo]{Lemma}

\newtheorem{Def}[theo]{Definition}
\newtheorem{prop}[theo]{Proposition}
\numberwithin{equation}{section}
\newtheorem{assump}{Assumption A-\hspace{-0.15cm}}
\theoremstyle{remark}
\newtheorem{rem}{Remark}[section]

\newcommand{\Zphi}{\bar{Z}^{\phi}_n}
\newcommand{\bdm}{\begin{displaymath}}
\newcommand{\edm}{\end{displaymath}}
\newcommand{\esp}{\mathbb{E}\,}
\newcommand{\N}{\mathbb{N}}
\newcommand{\Rd}{\mathbb{R}^d}
\newcommand{\la}{\langle}
\newcommand{\ra}{\rangle}
\newcommand{\epilim}{\operatornamewithlimits{epi-lim}}
\newcommand{\Dout}{D_{\infty,\mathrm{out}}}
\newcommand{\Din}{D_{\infty,\mathrm{in}}}
\newcommand{\D}{D_{\infty}}
\newcommand{\Cout}{C_{\infty,\mathrm{out}}}
\newcommand{\Cin}{C_{\infty,\mathrm{in}}}
\newcommand{\C}{C_{\infty}}
\newcommand{\ridom}{\mathrm{ri\,}\mathrm{dom\,}}

\newcommand{\X}{\mathcal X}

\begin{document}

\begin{abstract}
  We study in this article the large deviations for the weighted
  empirical mean $L_n =\frac{1}{n} \sum_1^{n} \mathbf{f}(x_i^n)\cdot
  Z_i ,$ where $(Z_i)_{i\in \N}$ is a sequence of $\Rd$-valued
  independent and identically distributed random variables with some
  exponential moments and where the deterministic weights
  $\mathbf{f}(x_i^n)$ are $m\times d$ matrices.  Here $\mathbf{f}$ is
  a continuous application defined on a locally compact metric space $({\mathcal X},
  \rho)$ and we assume that the empirical measure $
  \frac 1n \sum_{i=1}^n \delta_{x_i^n}$ weakly converges to some probability distribution $R$ with compact support ${\mathcal Y}$.\\
  The scope of this paper is to study the effect on the Large
  Deviation Principle (LDP) of \emph{outliers}, that is elements
  $x_{i(n)}^n \in \{x_i^n,\ 1\le i\le n\}$ such that
  $$ 
  \liminf_{n\rightarrow \infty} \rho(x_{i(n)}^n, {\mathcal Y})>0\ .
  $$ 
  We show that outliers can have a dramatic impact on the rate
  function driving the LDP for $L_n$.  We also
  show that the statement of a LDP in this case requires specific
  assumptions related to the large deviations of the single random
  variable $\frac{Z_1}n$. This is the main input with respect to a
  previous work by Najim \cite{Naj02}.
\end{abstract}
\maketitle
\bibliographystyle{plain}
\noindent {\sl Math. Subj. Class.:} Primary 60F10, Secondary 15A52, 15A18.

\section{Introduction}
\subsection*{The model}

\begin{figure}
\begin{center}
\epsfig{figure=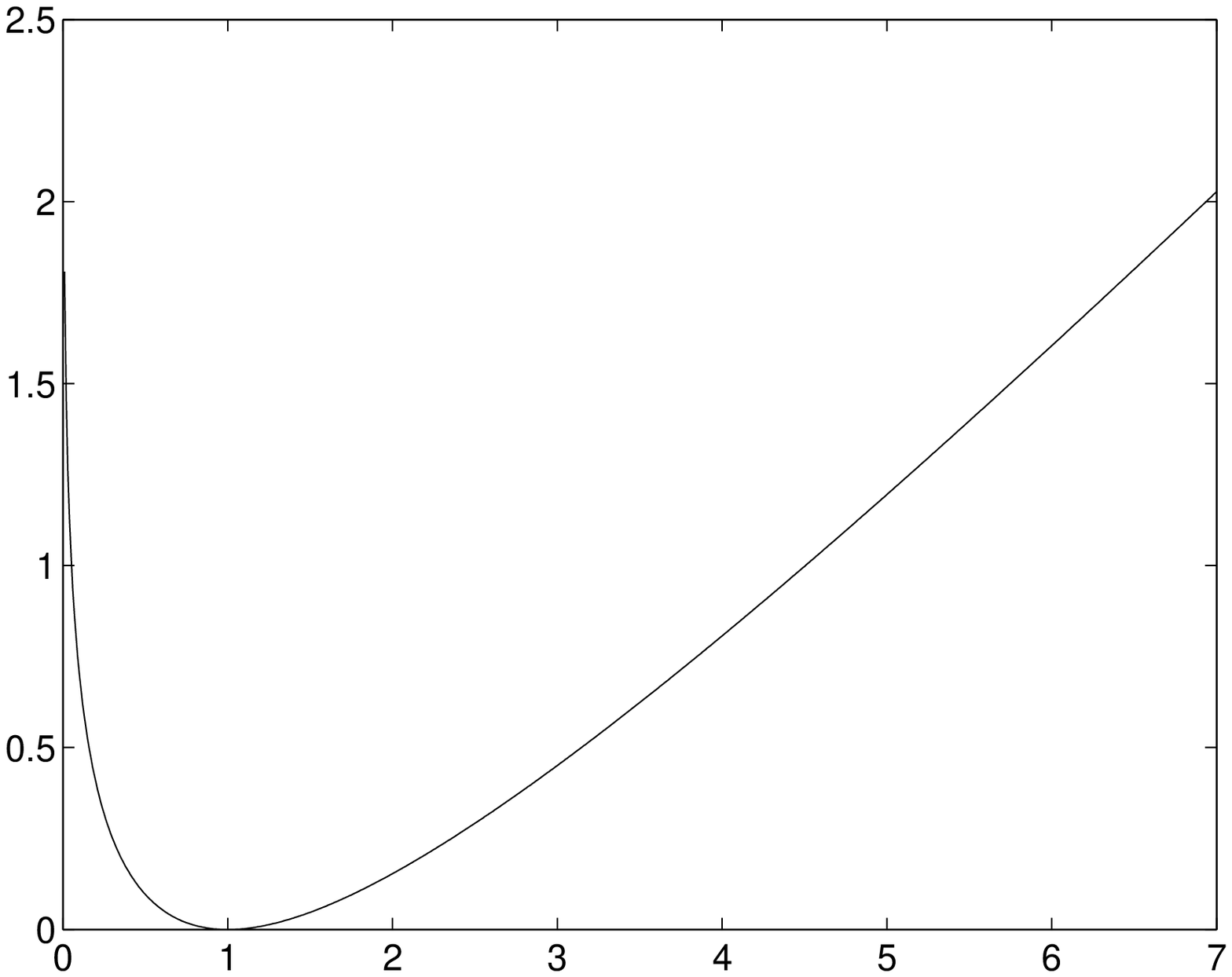,height=5cm}
\epsfig{figure=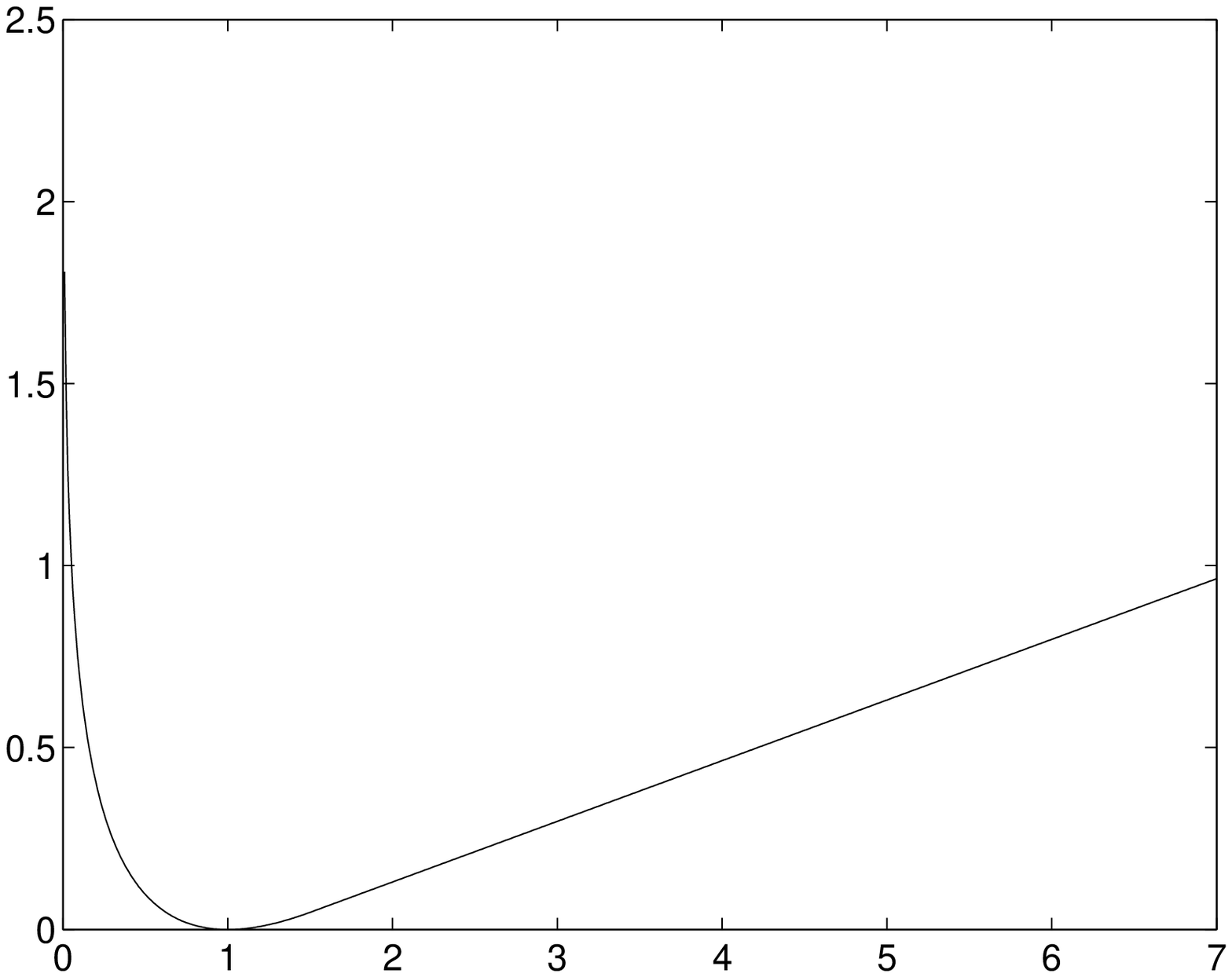,height=5cm}
\caption{\small The rate function of $\frac 1n \sum_{i=1}^n X_i^2$
  where the $X_i$'s are ${\mathcal N}(0,1)$ Gaussian i.i.d. random variables 
  (left); the rate function of $\frac 1n \sum_{i=1}^{n-1} X_i^2
  +\frac{3}{n} X_n^2$ (right). Both rate functions coincide for $x\le
  \frac32$ but the right one is linear for $x>\frac
  32$.}\label{dessin}
\end{center}
\end{figure}

We study in this article a Large Deviation Principle (LDP) for the weighted empirical mean
\begin{equation*}
L_n =\frac{1}{n}
\sum_1^{n} \mathbf{f}(x_i^n)\cdot Z_i ,
\end{equation*}
where $(Z_i)_{i\in \N}$ is a sequence of $\Rd$-valued independent and
identically distributed (i.i.d) random variables satisfying:
\begin{equation}\label{moment-exp}
\esp e^{\alpha |Z_1|}<\infty \qquad \textrm{for\ some}\quad \alpha>0.
\end{equation}
The application $\mathbf{f}:{\mathcal X}\rightarrow
\mathbb{R}^{m\times d}$ is a $m\times d$ matrix-valued continuous
function, $({\mathcal X},\rho)$ being a locally compact metric space. The term
$\mathbf{f}(x)\cdot Z$ denotes the product between matrix
$\mathbf{f}(x)$ and vector $Z$.  The set $\{x_i^n,1\le i\le n,\ n\ge
1\}$ is an ${\mathcal X}$-valued sequence of deterministic elements
such that the empirical measure
$\hat{R}_n\stackrel{\triangle}{=}\frac{1}{n} \sum_{i=1}^n
\delta_{x_i^n}$ satisfies:
\begin{equation}\label{mun}
\hat{R}_n  \xrightarrow[n\rightarrow\infty]{\mathrm{weakly}} R\ ,
\end{equation}
where $R$ is a probability measure with compact support ${\mathcal
  Y}$. \\
We focus in this paper on cases where there are outliers,
that is where some of the $x_i^n$ remain far from the support (also called bulk) of
$R$. Loosely speaking, one can think of an outlier as a sequence
$(x_{i(n)}^n, n\ge 1)$ satisfying:
\begin{equation}\label{loose-outlier}
\liminf_{n\rightarrow \infty} \rho(x_{i(n)}^n, {\mathcal Y})>0\ .
\end{equation}
At a large deviation level, such outliers may have a dramatic impact
on the shape of the rate function as demonstrated in the simple
example of {\sc Figure} \ref{dessin}. Although the model under study
looks very similar to the LDP studied in \cite{Naj02}, the presence of
outliers substantially modifies the resulting LDP and may naturally create
infinitely many non-exposed points (see the definition in \cite{DemZei98} and also Remarks \ref{non-expo} and \ref{non-expo-rem}) for the rate function.

The purpose of this article is to provide clear assumptions (which
cover situations where (\ref{loose-outlier}) can occur) over the set
$\{\mathbf{f}(x_i^n),\ 1\le i\le n,\ 1\le n\}$ and over $Z_i$ under
which fairly general LDP results can be proved.

\subsection*{Motivations and related work}
Such models are of particular interest in the field of statistical
mechanics (spherical spin glasses in \cite{BDG01}, spherical integrals
in the finite rank case in \cite{GuiMai05}, etc.)  where one has often
to establish a LDP for the empirical mean $L_n$ in the case where the
random variable $Z_i$ satisfies condition (\ref{moment-exp}). In
particular, spherical integrals are intimately connected to the study
of Deformed Ensembles (see \cite{Pec06} for instance for the
definition) in Random Matrix Theory.  In dimension one, $Z_i$ is
typically the square of a Gaussian random variable. The measure
$\frac{1}{n} \sum_{i=1}^n \delta_{x_i^n}$ is then a realization of the
empirical measure of the eigenvalues associated to a given random
matrix model and there are important cases when some of the $x_i^n$'s
stay far away from the support of $R$. Indeed, there has recently been
a strong interest in random matrix models (so-called spiked models)
where some of the largest eigenvalues lie out of the bulk, that is
where the set of limit points of $(x_i^n,\ 1\le i\le n,\ n\ge 1)$ can
differ from the support of $R$ (see Johnstone \cite{Joh01}, Baik et
al. \cite{BBP05}, \cite{BaiSil04pre}, P\'ech\'e \cite{Pec06}). These
spiked models are of particular interest for statistical applications
\cite{Joh01}.

The study of the LDP for weighted means was developed by Bercu et al.
\cite{BerGamRou97} for Gaussian functionals and considered in greater
generality in Najim \cite{Naj02}. In \cite{Naj02}, the LDP is stated
for $L_n$ under condition (\ref{moment-exp}) but in the case where
$(x_i^n,\ 1\le i\le n,\ 1\le n)$ is a subset of ${\mathcal Y}$, the
support of the limiting probability measure $R$. In particular, the
framework of \cite{Naj02} does not allow any of the $x_i^n$'s to lie
far from the bulk. LDPs involving outliers can be found in Bercu et
al. \cite{BerGamRou97},
Guionnet and Ma\"{\i}da \cite{GuiMai05}.
For related work concerning quadratic forms of Gaussian processes, we
shall also refer the reader to Bercu et al.  \cite{BerGamLav00},
Gamboa et al. \cite{GamRouZan99}, Bryc and Dembo \cite{BryDem97} and
Zani \cite{Zan99}.

\subsection*{Presentation of the results} The purpose of this article
is to establish the LDP for the empirical mean $L_n$ under the moment
assumption (\ref{moment-exp}) and under assumptions which allow the
presence of outliers (see (\ref{loose-outlier})). Such a LDP will rely
on the individual LDP for $\frac{Z_1}n$. This is the content of the
following assumption.
\begin{assump}\label{LDP-Particle}
The $\Rd$-valued random variable $Z_1$ satisfies the following exponential condition:
$$
\esp e^{\alpha |Z_1|}<\infty \qquad \textrm{for\ some}\quad \alpha>0,
$$
and 
$\frac {Z_1}n$ satisfies the LDP with a good rate function denoted by $I$.
\end{assump}
Note that if $\frac{Z_i}n$ does not satisfy a LDP, one can construct
counterexamples where $L_n$ does not fulfill a LDP (see for instance
\cite[Section 2.3]{Naj02}).  Finally, two subcases of Assumption
(A-\ref{LDP-Particle}) yield to two distinct classes of results:

\subsubsection*{The case where $I$ is convex (Assumption
  (A-\ref{convex-rf}), Section \ref{subsection-assump})} This paper is
mainly devoted to the study of this case.  If $I$ is convex then the
assumptions on the sets $C_n^{\mathbf{f}}=\{\mathbf{f}(x_i^n),\ 1\le
i\le n, \ 1\le n\}$ needed to state the LDP for $L_n$ are quite mild.
Apart from a standard compacity assumption (Assumption
(A-\ref{Compacity}), see Section \ref{subsection-assump}), the main
assumption over $C_n^{\mathbf{f}}$ (Assumption (A-\ref{Limit-Points}),
Section \ref{subsection-assump}) bears on the sole limiting points of
$C_n^{\mathbf{f}}$ (in the sense of Painlev\'e-Kuratowski convergence
of sets) and on their role in the LDP. It turns out that
(A-\ref{Limit-Points}) is an intricate assumption concerning the limiting
behaviour of $C_n^{\mathbf{f}}$ and some limiting points of
$C_n^{\mathbf{f}}$ involved in the definition of a certain convex
domain. This convex domain plays a role in the definition of the rate function of the LDP.
As demonstrated by examples in Section \ref{subsection-examples},
(A-\ref{Limit-Points}) covers a wide variety of models with outliers
in the convex case, at least those for which a LDP is to be expected.

Under Assumptions (A-\ref{LDP-Particle})-(A-\ref{Limit-Points}) and the
more classical assumption (A-\ref{hypo-EJP}) (convergence of $\hat
R_n$ to $R$), the empirical mean $L_n$ satisfies the LDP with a good
convex rate function (Theorem \ref{ldp-mem}). This rate function
admits a fairly good representation (in terms of convex features)
where the role of the outliers is quiet transparent (Theorem
\ref{theo-rf} and examples in Section \ref{section:examples}).

\subsubsection*{The case where $I$ is not convex} In this case, one can still prove 
the LDP but the assumptions over $C_n^{\mathbf{f}}$ are much more stringent and the rate function
is given by an abstract formula. Moreover, very few insight can be gained
by the study of the general formula of the rate function. It seems
that the study must be held on a case-by-case analysis.

\subsection*{Outline of the article}
In order to study the Large Deviations of $L_n$, we shall separate outliers from the bulk and split accordingly $L_n$ into two subsums:
\begin{eqnarray*}
L_n&=&
\frac 1n \sum_{\{x_i^n\ \textrm{far from the bulk}\}} \mathbf{f}(x_i^n)\cdot Z_i+
\frac 1n \sum_{\{x_i^n\ \textrm{near or in the bulk}\}} \mathbf{f}(x_i^n)\cdot Z_i\\
&\stackrel{\triangle}{=}& \pi_n + \tilde L_n.
\end{eqnarray*}
The idea is then to establish separately the LDP for each
subsum. 
This line of proof has been developed in the one-dimensional
setting for Gaussian quadratic forms by Bercu et al.
\cite{BerGamRou97} and is extended to the multidimensional
setting in this article. 

The paper is organized as follows. 
Sections \ref{section:partiel},
\ref{section:ldp} and \ref{section:examples} are devoted to the study
of the convex case. 

In Section \ref{section:partiel}, we study the Large Deviations for the following model:
\begin{eqnarray}\label{model-sec2}
\pi_n =\frac 1n \sum_{x_i^n \in C_n} \mathbf{f}(x_i^n)\cdot Z_i \quad \textrm{where}
\quad \frac{\mathrm{card}(C_n)}{n} \xrightarrow[n\rightarrow\infty]{} 0.
\end{eqnarray}
The main assumptions related to the set $C_n^{\mathbf{f}} =\{ \mathbf{f}(x_i^n);\ x_i^n \in C_n\}$
are stated and the LDP for $\pi_n$ is established.

In Section \ref{section:ldp}, the decomposition $L_n =\pi_n +\tilde
L_n$ where $\pi_n$ satisfies \eqref{model-sec2} is precisely specified, 
the LDP for $L_n$ is established and a
representation formula is given for the rate function.
Section \ref{section:examples} is devoted to examples of LDPs with
outliers in the convex case. 

A general LDP stated with an abstract rate function is established in
the non-convex case in Section \ref{section:nonconvex}.
In Section \ref{section:eigen}, a
partial study of the rate function is also carried out in the
non-convex case in the setting of a specific example.

Comments related to the link between the study of the spherical integral and 
the LDP of $L_n$ are made in Sections \ref{section:examples} (rank one
case) and \ref{section:eigen} (higher rank).

\section{The LDP for the partial mean $\pi_n$ in the convex case}\label{section:partiel}

Let $(C_n)_{n\geq 1}$ be a finite subset of ${\mathcal X}$. This section is devoted to the study of the LDP of 
$$
\pi_n = \frac 1n \sum_{x_i^n \in C_n} \mathbf{f}(x_i^n)\cdot Z_i \quad 
\textrm{where}\quad \frac{\mathrm{card}( C_n)}{n} \xrightarrow[n\rightarrow \infty]{} 0,
$$
with $\mathrm{card}(C_n)$ standing for the cardinality of the set $C_n$. It will be proved 
in Section \ref{decomposition} that $L_n$ can be decomposed as $\pi_n+\tilde L_n$ with $\pi_n$
as above.
\begin{rem} In the case where the random variable $Z_1$ satisfies
\begin{equation}\label{every-exp-moment}
\mathbb{E} e^{\alpha |Z_1|} <\infty\quad \textrm{for all} \quad \alpha \in \mathbb{R}^+,
\end{equation}
the following limit holds true:
$$
\limsup_{n\rightarrow \infty} \frac 1n\log \mathbb{P}\{ |\pi_n|>\delta \}=-\infty \quad \textrm{for all}\quad \delta>0.
$$
Otherwise stated $L_n$ and $\tilde L_n$ are exponentially equivalent and $\pi_n$ does not play any role 
at a large deviation level. Of course the situation is completely different if \eqref{every-exp-moment}
does not hold.
\end{rem}
We first introduce some notations as well as the concepts of inner limit, outer limit and
Painlev\'e-Kuratowski convergence for sets.
We then state the assumptions over the sets
$C_n^{\mathbf{f}}=\{ \mathbf{f}(x_i^n),\
x_i^n\in C_n\}$ and prove the LDP for $\pi_n$.

\subsection{Notations}
Denote by ${\mathcal B}({\mathcal Z})$ the Borel sigma-field of a given topological space ${\mathcal Z}$ 
(usually $\mathbb{R}^d$, $\mathbb{R}^m$, $\mathbb{R}^{m\times d}$ or ${\mathcal X}$). Denote by $|\cdot|$  a norm on 
any finite-dimensional vector space ($\mathbb{R}^d$, $\mathbb{R}^m$ or $\mathbb{R}^{m\times d}$).
In the sequel, we use bold letters $\mathbf{a},\mathbf{b},\mathbf{y}$, etc. to denote $m\times d$ matrices. We denote by 
$\la \cdot,\cdot\ra$  the scalar product in any finite-dimensional space and by $\cdot$
the product between vectors and  matrices with compatible size. 
Let $A$ be a subset of $\mathbb{R}^k$. We denote by $\bar{A}$ its closure, by $\mathrm{int}(A)$ its interior, 
by $\Delta(\cdot\mid A)$ the convex indicator function of the set $A$ and by $\Delta^*(\cdot\mid A)$ 
its convex conjugate (also called the support function of $A$), that is:
\begin{eqnarray*}
\Delta(\theta \mid A)&=&\left\{
\begin{array}{ll}
0&\textrm{if} \ \theta \in A,\\
\infty&\textrm{else}.
\end{array}\right. ,\\
\Delta^*(y \mid A)
&=& \sup_{ \theta \in \mathbb{R}^k}\{ \la y, \theta \ra -\Delta(\theta\mid A)\}
=\sup_{\theta\in A} \la y, \theta \ra,  
\end{eqnarray*}
where $y$ and $\theta$ are in $\mathbb{R}^k$. The following proposition whose proof is straightforward 
will be of constant use in the sequel.
\begin{prop}\label{support-function} 
Let $A$ be a subset of $\mathbb{R}^k$, then 
$$
\Delta^*(\cdot\mid A) = \Delta^*(\cdot \mid \bar{A}).
$$ 
If moreover $A$ is convex with non-empty interior, then 
$$
\Delta^*(\cdot\mid \mathrm{int}(A)) = \Delta^*(\cdot\mid A) = \Delta^*(\cdot \mid \bar{A}).
$$
\end{prop}

Let $D_n$ be a sequence of subsets of $\mathbb{R}^{m\times d}$. We
define its outer limit (denoted by $\Dout$) 
and its inner limit (denoted by $\Din$) 
by 
\begin{eqnarray*}
\Dout &=& \left\{
\mathbf{x}\in \mathbb{R}^{m\times d},\ \exists\, \phi:  \mathbb{N}
\rightarrow  \mathbb{N} \  \textrm{increasing,}\ 
\exists\, \mathbf{x}_{\phi(n)} \in D_{\phi(n)},\ \mathbf{x}_{\phi(n)}\xrightarrow[n\rightarrow \infty]{} \mathbf{x} \right\}\\  
\Din &=& \left\{
\mathbf{x}\in \mathbb{R}^{m\times d},\ \exists\, n_0,\ \forall\, n\ge n_0,\exists\, \mathbf{x}_n\in D_n,\ 
\mathbf{x}_{n}\xrightarrow[n\rightarrow \infty]{} \mathbf{x} \right\}\\
\end{eqnarray*}
The limit $\D$ of the sets $(D_n)$ exists if the outer limit and the inner limit are equal. 
Set convergence in this sense is known as Painlev\'e-Kuratowski convergence and in this case, we will denote:
$$
D_n\xrightarrow[n\rightarrow\infty]{\textrm{pk}} \D.
$$
For more details on Painlev\'e-Kuratowski convergence of sets, see Rockafellar and Wets \cite[Chapter 4]{RocWet98}.

\subsection{A preliminary analysis: Two simple examples}\label{subsection-examples} Consider
$$
C^{\mathbf{f}}_n=\left\{ \mathbf{f}(x_i^n),\ x_i^n \in C_n\right\}\quad \textrm{where}\quad
\frac{\mathrm{card}(C_n)}{n} \rightarrow 0.
$$
The sets $\Cin^{\mathbf{f}}$ and $\Cout^{\mathbf{f}}$ are respectively the inner and outer limits of $(C_n^{\mathbf{f}})$.
In the study of the forthcoming examples, we will focus on the links between the LDP for
$\pi_n$ and the sets $\Cin^{\mathbf{f}}$ and $\Cout^{\mathbf{f}}$. This section is aimed at introducing Assumption 
(A-\ref{Limit-Points}) but can be skipped as no further notation is introduced.

\subsubsection{Example 1: A simple case where the LDP fails to hold for $\pi_n$}\label{example-1}
Let $X$ be a standard Gaussian random variable and consider
$\pi_n=\frac{2+(-1)^n}{n} X^2$.  Direct computations yield the LDP for
$\pi_{2n}$ (resp. $\pi_{2n+1}$) with good rate function
$\Delta^*_{\mathrm{even}}$ (resp.  $\Delta^*_{\mathrm{odd}}$) where
$$
\Delta^*_{\mathrm{even}}(z)=\left\{
\begin{array}{ll}
z/6 & \mathrm{if}\ z>0,\\
\infty & \mathrm{else}.
\end{array}\right. 
\quad \textrm{and} \quad
\Delta^*_{\mathrm{odd}}(z)=\left\{
\begin{array}{ll}
z/2 & \mathrm{if}\ z>0,\\
\infty & \mathrm{else}.
\end{array}\right. 
$$
Therefore one cannot expect the LDP for $(\pi_n,n\in \N)$.

\subsubsection{Example 2:  The LDP holds after modification of Example 1}\label{example-2}
Let $X$ and $Y$ be independent standard Gaussian random variables and consider 
$\pi_n=\frac{2+(-1)^n}{n} X^2 +\frac 4n Y^2$. 
In this case, $\pi_{2n}$ and 
$\pi_{2n+1}$ satisfy the LDP (by a direct analysis) with the same rate function 
$$
\Delta^*(z)=\left\{
\begin{array}{ll}
z/8 & \mathrm{if}\ z>0,\\
\infty & \mathrm{else}.
\end{array}\right.
$$
This yields the LDP for the whole sequence $(\pi_n,n\in\N)$ with rate function $\Delta^*$. 

Despite the erratic behaviour of $\frac{2+(-1)^n}{n} X^2$ (as seen in
the previous example), the LDP holds due to presence of the term $\frac 4n Y^2$.

\subsubsection{Comparison of the two examples}
Denote by 
$$
{\mathcal D}_y=\{\lambda\in \mathbb{R},\ \log \mathbb{E} e^{\lambda y X^2} <\infty\}
=\left(-\infty,(2y)^{-1}\right)
$$ 
where $X$ is a standard 
Gaussian random variable.

In the case of Example 1, one can easily check that $C_{2n}^{\mathbf{f}}= \{3\}$ and 
$C_{2n+1}^{\mathbf{f}}= \{1\}$. Thus $\Cout^{\mathbf{f}}=\{1,3\}$ while
$\Cin^{\mathbf{f}}=\emptyset$. It is straightforward to check that the rate functions 
driving the LDP of $\pi_{2n}$ and $\pi_{2n+1}$ can be expressed as:
$$
\Delta^*_{\mathrm{even}}(z)\ =\ \sup_{\lambda\in {\mathcal D}_3} \lambda z\qquad\mathrm{and}\qquad 
\Delta^*_{\mathrm{odd}}(z)\ =\ \sup_{\lambda\in {\mathcal D}_1} \lambda z,
$$ 
The very reason for which the LDP does not hold in this case is that
$$
\bigcap_{y\in \Cout^{\mathbf{f}}} {\mathcal D}_y\  \neq\  
\bigcap_{y\in \Cin^{\mathbf{f}}} {\mathcal D}_y.
$$

In the case of Example 2, $C_{2n}^{\mathbf{f}}=\{3,4\}$ while 
$C_{2n+1}^{\mathbf{f}}=\{1,4\}$. Therefore $\Cout^{\mathbf{f}}=\{1,3,4\}$ while 
$\Cin^{\mathbf{f}}=\{4\}$. Despite the fact that $\Cout^{\mathbf{f}} \neq \Cin^{\mathbf{f}}$,
the LDP holds in this case with good rate function given by:
\begin{eqnarray*}
\Delta^*(z)&=& \sup_{\lambda\in {\mathcal D}_4} \lambda z.
\end{eqnarray*} 
As we shall see, the underlying reason for which the LDP holds is
$$
\bigcap_{y\in \Cout^{\mathbf{f}}} {\mathcal D}_y\  =\  
\bigcap_{y\in \Cin^{\mathbf{f}}} {\mathcal D}_y\ \left(=\ {\mathcal D}_4\right),
$$
and this will be a key-point in the statement of Assumption (A-\ref{Limit-Points}).

We are now in position to state the assumptions and the main result. 
\subsection{Assumptions and main results}\label{subsection-assump}
Let $C_n$ be a finite subset of ${\mathcal X}$ and recall that
$$
C^{\mathbf{f}}_n=\left\{ \mathbf{f}(x_i^n),\ x_i^n \in C_n \right\}\quad \textrm{where}\quad \frac{\mathrm{card}(C_n)}{n}
\xrightarrow[n\rightarrow \infty]{}0.
$$
Let $\mathbf{y}$ be a
$m\times d$ matrix and denote by
\begin{equation}\label{def-D}
{\mathcal D}_{\mathbf{y}}=\left\{ \lambda \in \mathbb{R}^m,\ \log \esp e^{\la \lambda,  \mathbf{y}\cdot Z_1\ra} <\infty \right\}.
\end{equation}
We can now state our assumptions. 

Assume that $Z_1$ is a $\mathbb{R}^d$-valued random variable satisfying Assumption (A-\ref{LDP-Particle})
and recall that $I$ is the rate function associated to $\frac {Z_1} n$.

\begin{assump}\label{convex-rf}
Let 
$ \mathcal D_Z \stackrel{\triangle}{=} \{\theta \in
\Rd,\  \log \esp e^{\langle \theta, Z_1\rangle}  < \infty\},
$
then 
$$
I(z) = \Delta^*(z \mid \mathcal D_Z).
$$
In particular, $I$ is a convex rate function.
\end{assump}

\begin{assump}\label{Compacity} 
Let $(D_n)_{n\ge 1}$ be a sequence of non empty subsets of
  $\mathbb{R}^{m\times d}$. There exists a compact set $K\subset \mathbb{R}^{m\times d}$ such that
$D_n \subset K$ for every $n\ge 1$.
\end{assump}
\begin{rem}
This assumption implies  in particular that 
the outer limit $\Dout$ of $(D_n)_{n\ge 1}$ is a nonempty compact set of $\mathbb{R}^{m\times d}$.
\end{rem}

\begin{assump}\label{Limit-Points}
Let  $(D_n)_{n\ge 1}$ be a sequence of subsets of
  $\mathbb{R}^{m\times d}$. Denote by $\Din$ and $\Dout$ its inner
  and outer limits. Then:
$$
\bigcap_{\mathbf{y}\in \Din}{\mathcal D}_{\mathbf{y}}
=\bigcap_{\mathbf{y}\in \Dout}{\mathcal D}_{\mathbf{y}}
$$
where ${\mathcal D}_{\mathbf{y}}$ is defined by (\ref{def-D}). 
\end{assump}

\begin{rem}
If $(D_n)_{n\ge 1}$ fulfills (A-\ref{Compacity}) and
(A-\ref{Limit-Points}), then in particular, $\Din$ is not empty.
\end{rem}
We can now state the main result of the section.
\begin{theo}\label{ldp}
  Assume that $(Z_i)_{i\in \mathbb{N}}$ is a sequence of
  $\mathbb{R}^d$-valued i.i.d. random variables. Assume moreover that
  (A-\ref{LDP-Particle}) and (A-\ref{convex-rf}) hold for $Z_1$.
  Assume that $({\mathcal X},\rho)$ is a metric space and let
    $C_n\subset{\mathcal X}$ be such that
$$
\frac{\mathrm{card}(C_n)}{n}
\xrightarrow[n\rightarrow \infty]{}0.
$$
Denote by $C_n^{\mathbf{f}}=\{ \mathbf{f}(x_i^n),\ x_i^n \in
  C_n\}$ where $\mathbf{f}:{\mathcal X} \rightarrow \mathbb{R}^{m\times d}$ is continuous. 
  Assume that (A-\ref{Compacity}) and (A-\ref{Limit-Points})
  hold for the sequence of sets $(C^{\mathbf{f}}_n)_{n\in
    \mathbb{N}}$. Then the random variable
$$
\pi_n =\frac{1}{n} \sum_{x_i^n \in C_n} \mathbf{f}(x_i^n)\cdot Z_i
$$
satisfies the LDP in $(\mathbb{R}^m,{\mathcal B}(\mathbb{R}^m))$ with good rate function 
$$
\Delta^*(z\mid {\mathcal D})=\sup_{} \{\la \lambda ,z \ra,\ \lambda \in {\mathcal D}\}
\qquad\textrm{where}\qquad 
{\mathcal D}=\bigcap_{\mathbf{y}\in \Cin^{\mathbf{f}}}{\mathcal D}_{\mathbf{y}}
=\bigcap_{\mathbf{y}\in \Cout^{\mathbf{f}}}{\mathcal D}_{\mathbf{y}}.
$$
\end{theo}

\begin{rem}[On Assumption (A-\ref{Limit-Points})]
  A close look to the proof of Theorem \ref{ldp} shows
  that the rate function that drives the lower bound of the LDP is the support function of 
  $\cap_{\mathbf{y}\in
      \Cin^{\mathbf{f}}}{\mathcal D}_{\mathbf{y}}$ while the
  rate function that drives the upper bound is the support function of $ \cap_{\mathbf{y}\in \Cout^{\mathbf{f}}}{\mathcal
      D}_{\mathbf{y}}$. Both rate functions coincide when assuming (A-\ref{Limit-Points}).
(see also the examples in Section \ref{subsection-examples}).
\end{rem}



\subsection{Proof of Theorem \ref{ldp}}

In order to prove Theorem \ref{ldp} , we follow the strategy developed in \cite{Naj02}, essentially based on 
an exponential approximation technique. The next proposition is the counterpart of Lemma 5.1 in \cite{Naj02}.

\begin{lemma}
\label{nw}
Let $\phi : \mathbb{N} \setminus\{0\} \rightarrow \mathbb{N} \setminus \{0\}$ be such that $\frac{\phi(n)}{n} \xrightarrow[n\rightarrow \infty]{}
  0$. Let $(Z_i)$ be a sequence of $\Rd$-valued random variables
  satisfying (A-\ref{LDP-Particle}) and (A-\ref{convex-rf}). Then 
$\Zphi\stackrel{\triangle}{=}\frac 1n \sum_{i=1}^{\phi(n)} Z_i$ satisfies the LDP in $\Rd$
with good rate function
 given by 
$$ I(y)=\Delta^*(y\mid {\mathcal D}_Z)$$
where ${\mathcal D}_Z$ is defined in (A-\ref{convex-rf}).
\end{lemma}

\begin{proof}
Denote by $\Lambda^{\phi}_n$ the log-Laplace transform of $\Zphi$, 
i.e. $\Lambda^{\phi}_n(\theta)=\log \esp e^{\la \theta,\bar{Z}^{\phi}_n\ra}$. Then 
$$
\frac{1}{n} \Lambda^{\phi}_n(n\theta)=\frac{\phi(n)}{n}\log \esp\, e^{\la \theta, Z_i\ra}\xrightarrow[n\rightarrow \infty]{} 
\Delta(\theta\mid{\mathcal D}_Z).
$$
Therefore, the large deviation upper bound holds for $\Zphi$ with rate function $I$ 
by Theorem 2.3.6 (a) in \cite{DemZei98}. 
To prove the large deviation lower bound, it is sufficient to prove that
$$
-I(y) \le \liminf_{n\rightarrow \infty} \frac{1}{n}\log
\mathbb{P}\left(\Zphi\in B(y,\varepsilon)\right)
$$
where $B(y,\varepsilon)=\{y'\in \Rd,\ |y'-y|<\varepsilon\}$.
Define
$$
\tilde{Z}_n^{\phi}=
\left\{
\begin{array}{ll}
\frac 1n \sum_{i=2}^{\phi(n)} Z_i &  \textrm{if}\ \phi(n)\ge2,\\
0 & \textrm{otherwise}.
\end{array}
\right. . 
$$
Then $
\{Z_1/n \in B(y,\varepsilon/3)\} \cap \{\tilde{Z}^{\phi}_{n}\in B(0,\varepsilon/3)\}
\subset \{\Zphi\in B(y,\varepsilon)\}
$
which yields 
\begin{multline}\label{eq:min}
\frac{1}{n}\log\mathbb{P} \left(Z_1/n \in B(y,\varepsilon/3)\right)
+  \frac{1}{n}\log \mathbb{P} \left(\tilde{Z}^{\phi}_{n}\in B(0,\varepsilon/3)\right)\\
\le \frac{1}{n}\log \mathbb{P}\left(\Zphi\in B(y,\varepsilon)\right).
\end{multline}
Exponential Markov inequality yields $\lim_{n\rightarrow \infty}
\mathbb{P}\{|\tilde{Z}^{\phi}_{n}|>\varepsilon/3\}=0$ which readily
implies that $\lim_{n\rightarrow \infty}\mathbb{P} \{\tilde{Z}^{\phi}_{n}\in B(0,\varepsilon/3)\}=1.$
Consequently, taking the liminf in both sides of
(\ref{eq:min}) and using the lower bound for the single variable
$\frac{Z_1}{n}$ yields the desired lower bound. The proof is
completed. 
\end{proof}
We first consider Theorem \ref{ldp} under an additional assumption.
\begin{lemma}
\label{ldp-pk}
Under the same assumptions as in Theorem \ref{ldp} and if we assume in addition
that
\begin{equation}\label{tight-assumption}
C_n^{\mathbf{f}}
  \xrightarrow[n\rightarrow\infty]{\mathrm{pk}}
  C^{\mathbf{f}}_{\infty},
\end{equation}
then $\pi_n$ satisfies the LDP in $\Rd$ with good rate function
$\Delta^*(\,\cdot \mid \mathcal D),$
where ${\mathcal D}=\cap_{\mathbf{y}\in C_{\infty}^{\mathbf{f}}}{\mathcal D}_{\mathbf{y}}$.
\end{lemma} 

\noindent Proof of Lemma \ref{ldp-pk} is postponed to Appendix \ref{proof-ldp-pk}.\\

\noindent We now relax the extra assumption (\ref{tight-assumption}) and prove
Theorem \ref{ldp}.  The scheme of the proof is the following. We first
show, using directly the result in Lemma \ref{ldp-pk}, that the lower
bound is driven by the support function of the set $\bigcap_{\mathbf y
  \in \Cin^{\mathbf f}} \mathcal D_{\mathbf y}$. We then obtain that the
upper bound is driven by the support function of the set
$\bigcap_{\mathbf y \in \Cout^{\mathbf f}} \mathcal D_{\mathbf y}$, by
majorizing the log-Laplace of $\pi_n$.  Under Assumption
(A-\ref{Limit-Points}), both bounds coincide and we get the full LDP.

\begin{proof}[Proof of Theorem \ref{ldp}]
To get the lower bound, we split $C_n^{\mathbf f}$ into two disjoint subsets:
\begin{equation}\label{split}
C_n^{\mathbf{f}}={\mathcal I}_n^{\mathbf{f}} \cup {\mathcal O}_n^{\mathbf{f}}\qquad
\textrm{ where }\qquad {\mathcal I}_n^{\mathbf{f}}\xrightarrow[n\rightarrow\infty]{\mathrm{pk}} \Cin^{\mathbf{f}}
\end{equation}
Let us sketch the construction of  ${\mathcal
  I}_n^{\mathbf{f}}$. Let $B(z,\frac 1m)$ be a ball centered in $z\in \Cin^{\mathbf{f}}$ with radius $\frac 1m$. 
Since $\Cin^{\mathbf{f}}$ is compact by (A-\ref{Compacity}), there exist $\left(z_{\ell}\right)_{1\le \ell\le L_m}$ such that
$$
\Cin^{\mathbf{f}} \subset \bigcup_{\ell =1}^{L_m}B\left(z_{\ell},\frac 1m\right)\quad \text{ and }\quad B\left(z_{\ell},\frac 1m\right)\cap \Cin^{\mathbf{f}}\not= \emptyset\quad \textrm{for}\quad 1\le \ell\le L_m. 
$$ 
The mere definition of $\Cin^{\mathbf{f}}$ yields that there exists $\psi(m)$
such that for all $\ell,\ 1\le \ell\le L_m$:
$$
\forall n\ge \psi(m),\quad  \exists \mathbf{f}(x_{i_{\ell}}^n)\in 
B\left(z_{\ell},\frac 1m\right)
\quad \textrm{with}\quad \mathbf{f}(x_{i_{\ell}}^n)\in C_n^{\mathbf{f}}.
$$
Denote by ${\mathcal A}_{n,m}$ ($n\ge \psi(m)$) such a collection of
$\mathbf{f}(x_{i_{\ell}}^n)$'s. Choose now similarly a collection of
balls with radius $\frac 1{m+1}$ and the related $\psi(m+1)$ with
$\psi(m+1)>\psi(m)$,
and set
$$
{\mathcal I}_n^{\mathbf{f}}={\mathcal A}_{n,m}\quad \textrm{if}\quad
\psi(m)\le n< \psi(m+1).
$$
With such a definition, it is straightforward to check that ${\mathcal I}_n^{\mathbf{f}}\xrightarrow[]{\mathrm{pk}} \Cin^{\mathbf{f}}.$
We write 
\begin{eqnarray*}
\pi_n &=&\frac 1n \sum_{x_i^n \in \mathbf{f}^{-1}({\mathcal
  I}_n^{\mathbf{f}})} \mathbf{f}(x_i^n) \cdot Z_i 
+ \frac 1n \sum_{x_i^n \notin \mathbf{f}^{-1}({\mathcal
  I}_n^{\mathbf{f}})} \mathbf{f}(x_i^n) \cdot Z_i\ , \\
&\stackrel{\triangle}{=}& \pi_n^{\mathcal
  I} + \pi_n^{\mathcal O}\ .
\end{eqnarray*}
The lower bound can be established as in Lemma \ref{nw}.
Let us prove that:
\begin{equation}\label{borneinf}
-\Delta^*(z\mid \cap_{\mathbf y \in \Cin^{\mathbf f}} \mathcal
 D_{\mathbf y}) \le \liminf_{n\rightarrow \infty} \frac 1n \log
 \mathbb{P}  \left(\pi_n \in B(z,\varepsilon)\right).
\end{equation}
Since 
$$
\{ \pi_n^{\mathcal I}\in B(z,\varepsilon/3) \} \cap 
\{ \pi_n^{\mathcal O}\in B(0,\varepsilon/3) \}\subset \{ \pi_n\in B(z,\varepsilon) \},
$$
one has 
\begin{multline}\label{minoration1}
\frac{1}{n}\log\mathbb{P}  \left( \pi_n^{\mathcal I} \in B(z,\varepsilon/3)\right)
+  \frac{1}{n}\log \mathbb{P}  \left(\pi_n^{\mathcal O}\in B(0,\varepsilon/3)\right)\\
\le \frac{1}{n}\log \mathbb{P} \left(\pi_n \in B(z,\varepsilon)\right).
\end{multline}
Exponential Markov inequality yields $\lim_{n\rightarrow \infty} \mathbb{P}(|\pi_n^{\mathcal O}|>\varepsilon/3)=0$.
This in turn implies that $\lim_{n\rightarrow \infty}\mathbb{P} \left(\pi_n^{\mathcal O}\in B(0,\varepsilon/3)\right)=1$.
Since $\pi_n^{\mathcal I}$ fulfills assumptions of Lemma \ref{ldp-pk}, the following lower bound holds:
\begin{equation}\label{minoration2}
 -\Delta^*\left(z\mid \cap_{\mathbf y \in \Cin^{\mathbf f}} \mathcal D_{\mathbf y} \right) \le \frac{1}{n}\log\mathbb{P} \left( \pi_n^{\mathcal I} \in B(z,\varepsilon/3)\right) 
\end{equation}
Consequently, taking the liminf in both sides of (\ref{minoration1}) and using (\ref{minoration2})
yields the desired lower bound. The proof of the lower bound is completed.

Let us now prove the upper bound.
Denote by $\Lambda_n(\lambda)$ the log-Laplace transform of $\pi_n$, i.e. $\Lambda_n(\lambda)=\log \esp e^{\la \lambda, \pi_n \ra}$.
In order to prove the upper bound, we estimate the following limit:
$$
\frac 1n \Lambda_n(n \lambda)= \frac 1n \sum_{x_i^n\in C_n} \log \esp e^{\la \lambda , \mathbf{f}(x_i^n) \cdot Z_i\rangle}
\qquad \textrm{where}\qquad \frac{\mathrm{card}(C_n)}n\xrightarrow[n\rightarrow \infty]{}0.
$$
We shall prove that 
\begin{equation}\label{ub}
\limsup_{n\rightarrow \infty} \frac 1n \Lambda_n( \lambda) \le
\Delta(\lambda\mid \mathrm{int}(\cap_{\mathbf y \in \Cout^{\mathbf f}} \mathcal D_{\mathbf y})).
\end{equation}
Theorem 4.5.3 in \cite{DemZei98} will then yield:
\begin{eqnarray}
\limsup_{n\rightarrow \infty} \frac 1n \log \mathbb P(\pi_n \in F) & \le & - \inf_{z
  \in F} \Delta^*(z\mid \mathrm{int}(\cap_{\mathbf y \in \Cout^{\mathbf f}} \mathcal D_{\mathbf y}))\nonumber \\
& \stackrel{(a)}{=}& - \inf_{z
  \in F} \Delta^*(z\mid \cap_{\mathbf y \in \Cout^{\mathbf f}} \mathcal D_{\mathbf y}) \label{ldp-upperbound}
\end{eqnarray}
for any closed set $F$. Equality $(a)$ follows from  Proposition
\ref{support-function} and the fact that
$\mathrm{int}(\cap_{\mathbf y \in \Cout^{\mathbf f}} \mathcal
D_{\mathbf y})$ is a non-empty convex set due to (A-\ref{LDP-Particle}).

In order to prove \eqref{ub}, consider $\lambda \in \Rd$ such that
\begin{equation}\label{limsup-strict-pos} 
\limsup_{n\rightarrow \infty} \frac 1n \Lambda_n(n \lambda) >0.
\end{equation}
>From \eqref{limsup-strict-pos}, we can successively:
\begin{itemize}
\item[-] extract a subsequence $n_{\alpha}$ from $n$ such that 
$$
\lim_{n \rightarrow \infty} \frac 1{n_{\alpha}}
\sum_{x_i^{n_{\alpha}}\in C_{n_{\alpha}}} \log \mathbb{E} e^{\langle
  \lambda, \mathbf{f}(x_i^{n_{\alpha}}) \cdot Z_i \rangle} >0;
$$
\item[-] extract a subsequence $n_{\beta}$ from $n_{\alpha}$ such that 
$$
\lim_{n \rightarrow \infty} \mathbb{E}e^{\langle \lambda,
  \mathbf{f}(x_i^{n_{\beta}}) \cdot Z_i \rangle} =\infty,
$$
\item[-] extract a subsequence $n_{\gamma}$ from $n_{\beta}$ such that
$$
\mathbf{f}(x_i^{n_{\gamma}}) \xrightarrow[n \rightarrow \infty]{} \mathbf{y}_0.
$$
One can notice in particular that $\mathbf{y}_0\in \Cout^{\mathbf{f}}$.
\end{itemize}
Let us now prove that 
\begin{equation}\label{lambda-frontier}
\lambda\notin \mathrm{int}({\mathcal D}_{\mathbf{y}_0}).
\end{equation}
Assume that (\ref{lambda-frontier}) is not true.
Then there exists $p>1$ such that $p\lambda \in {\mathcal D}_{\mathbf{y}_0}$. Let $\varepsilon>0$ 
be arbitrarily small. Then, if $n$ is large enough to ensure that 
$|\lambda| |\mathbf{f}(x_i^{n_{\gamma}}) -\mathbf{y}_0| \le \varepsilon/q$ where $1/p +1/q=1$, one has 
\begin{eqnarray*}
\mathbb{E}\,e^{\langle \lambda, \mathbf{f}(x_i^{n_{\gamma}}) \cdot Z\rangle }
&=&\mathbb{E}\,e^{\langle \lambda, \mathbf{y}_0 \cdot Z\rangle } 
e^{\langle \lambda, (\mathbf{f}(x_i^{n_{\gamma}})-\mathbf{y}_0) \cdot Z\rangle }\\
&\le& \left( \mathbb{E}\,e^{p \langle \lambda, \mathbf{y}_0 \cdot Z\rangle } \right)^{\frac 1p}
\left( \mathbb{E}\,e^{\varepsilon |Z| } \right)^{\frac 1q}.
\end{eqnarray*}
This contradicts the fact that 
$$
\lim_{n\rightarrow \infty} \mathbb{E}e^{\langle \lambda,  \mathbf{f}(x_i^{n_{\gamma}}) Z_i \rangle} =\infty.
$$ 
\noindent Therefore \eqref{lambda-frontier} holds and yields  that
$\lambda\notin \mathrm{int}(\cap_{\mathbf{y}\in
  \Cout^{\mathbf{f}}}{\mathcal D}_{\mathbf{y}_0})$. From this, we
deduce that
$$
\limsup_{n\rightarrow \infty} \frac 1n \Lambda_n(n \lambda) >0 \quad \Rightarrow \quad \lambda\notin 
\mathrm{int}(\cap_{\mathbf{y}\in \Cout^{\mathbf{f}}}{\mathcal D}_{\mathbf{y}}).
$$
Otherwise stated:
$$
\limsup_{n\rightarrow \infty} \frac 1n \Lambda_n(n \lambda) \le
\Delta\left(\lambda \mid \mathrm{int}(\cap_{\mathbf{y}\in
  \Cout^{\mathbf{f}}}{\mathcal D}_{\mathbf{y}})\right).
$$
Therefore, \eqref{ub} is proved and so is \eqref{ldp-upperbound}. 

Gathering the lower bound \eqref{borneinf}, the upper bound
\eqref{ldp-upperbound} and Assumption (A-\ref{Limit-Points}) yield the
full LDP for $\pi_n$.
\end{proof}

\section{The LDP for the empirical mean and the rate function in the convex case}\label{section:ldp}

Our goal is now to get the full LDP for $L_n$ (Theorem \ref{ldp-mem}
below). As announced in the outline of the article, the first step is
to split the $x_i^n$'s into two different subsets according to whether
they live near the support of the limiting measure or whether they are outliers.

\subsection{The decomposition ${\bf L_n=\pi_n + \tilde L_n}$}\label{decomposition}
Recall that $({\mathcal X},\rho)$ is a metric space.
\begin{prop}\label{proper-decomposition} Let  $A_n=\{x_i^n,\ 1\le i\le n\}$.
Assume that 
$$
\hat R_n =\frac 1n \sum_{i=1}^n \delta_{x_i^n} \xrightarrow[n\rightarrow\infty]{\mathrm{weakly}} R.
$$
and denote by ${\mathcal Y}$ the support of $R$. Then there exist 
subsets $B_n$ and $C_n=A_n \setminus B_n$ such that
\begin{enumerate}
\item $\frac{\mathrm{card}(B_n)}{n} \xrightarrow[n\rightarrow \infty]{} 1$,
  \label{prop1}
\item $\frac{1}{\mathrm{card}(B_n)} \sum_{x_i^n \in B_n} \delta_{x_i^n}
  \xrightarrow[n\rightarrow \infty]{\mathrm{weakly}} R$, \label{prop2}
\item $ \rho(B_n, {\mathcal Y}) \xrightarrow[n\rightarrow \infty]{} 0$
  where ${\mathcal Y}$ is the support of $R$. \label{prop3}
\end{enumerate} 
\end{prop}

We will then set 
$$
\tilde L_n =\frac 1n \sum_{x_i^n \in B_n} \mathbf{f}(x_i^n)\cdot Z_i \quad 
\textrm{and}\quad \pi_n =\frac 1n \sum_{x_i^n \in C_n} \mathbf{f}(x_i^n)\cdot Z_i.
$$
Note that since $\mathrm{card}(B_n)+\mathrm{card}(C_n)=n$, property
\eqref{prop1} yields then that
$
\frac{\mathrm{card}(C_n)}{n} \rightarrow 0
$
as $n$ goes to infinity.
\begin{proof}
{\em Construction of $B_n$}. Let $m\ge 1$ be fixed and
denote by ${\mathcal Y}_m$ the $\frac 1m$-blowup of ${\mathcal Y}$,
i.e.  ${\mathcal Y}_m=\{ x\in \X,\ \rho(x,{\mathcal Y})<\frac 1m \}$
where ${\mathcal Y}$ is the support of $R$. Then $\frac 1n \sum_1^n
1_{{\mathcal Y}_m}(x_i^n) \rightarrow 1$; in particular there exists
$\psi_m\ge 1$ such that for all $n\ge \psi_m$:
$$
\left| \frac 1n \sum_{i=1}^n 1_{{\mathcal Y}_m} (x_i^n) -1\right| <\frac 1m.
$$
One can then build recursively a sequence of integers $(\psi_m)_{m\in
  \mathbb{N}}$ such that $\psi_m<\psi_{m+1}$ (so that
$\psi_m\rightarrow \infty$ as $m\rightarrow \infty$). Set
$$
B_n=\{x_i^n \in {\mathcal Y}_m,\ 1\le i\le n\} \quad \textrm{for}\quad \psi_m\le n<\psi_{m+1}.
$$ 
We prove property (\ref{prop1}) and leave the proofs of
properties (\ref{prop2}) and (\ref{prop3}) to the reader.  

Let $\varepsilon>0$ be fixed and take $m$ such that $\frac 1m < \varepsilon$.
For such an $m$, take the corresponding $\psi_m$ and let $n\ge \psi_m$. Then,
$$
\left| \frac{\mathrm{card}(B_n)}{n} -1\right| =\left|
  \frac{\sum_{i=1}^n 1_{{\mathcal Y}_m}(x_i^n)}{n} -1\right| \le \frac
1m <\varepsilon.
$$
Since $\varepsilon>0$ is arbitrary, property (\ref{prop1}) is proved. 
\end{proof}

\subsection{The LDP for the empirical mean $L_n$}
In order to get the full LDP for $L_n =\tilde L_n +\pi_n$, we need to
prove the LDP for $\tilde L_n$. We will mainly rely on the results in
\cite{Naj02}. The following assumption is needed:
\begin{assump}\label{hypo-EJP}
Assume that $({\mathcal X}, \rho)$ is a locally compact metric space.
The family $(x_i^n,1\le i\le n,n\ge 1)\subset {\mathcal X}$ satisfies 
$$
\hat R_n = \frac 1n \sum_{i=1}^n \delta_{x_i^n} \xrightarrow[n\rightarrow\infty]{\textrm{weakly}} R,
$$
where $R$ is a probability measure over $({\mathcal X}, {\mathcal
  B}({\mathcal X}))$. Moreover, the support of $R$
denoted by ${\mathcal Y}$ is a compact set and for every non-empty open set $U$ of ${\mathcal Y}$ (for the 
induced topology over ${\mathcal Y}$), $R(U)>0$.
\end{assump}
\begin{rem} The LDP may fail to hold if the last part of Assumption
  (A-\ref{hypo-EJP}), that is $R(U)>0$ for $U$ non-empty
  open set, is not fulfilled. Counterexamples, also closely related to Assumption
  (A-\ref{LDP-Particle}), are developed in \cite{Naj02}.
\end{rem}

We recall that we denote by $\Lambda(\theta)=\log \esp e^{\la \theta,
  Z_1\ra }$ the log-Laplace transform of $Z_1$. We introduce the
following functional
\begin{equation} \label{Gamma}
\Gamma(\lambda)=\int_{\mathcal X} \Lambda\left(\sum_{k=1}^m \lambda_k f_k(x)\right) R(dx),
\end{equation}
where $\lambda=(\lambda_1,\cdots,\lambda_m)\in \mathbb{R}^m$ and $f_k$
denotes the k$^{\textrm{th}}$ row of matrix $\mathbf{f}$.  Let
$\Gamma^*$ be the convex conjugate of $\Gamma$:
$$
\Gamma^*(z)=\sup_{\lambda\in \mathbb{R}^m}\left\{ \la \lambda, z \ra- \Gamma(\lambda)\right\}. 
$$

We can now state the LDP.

\begin{theo}\label{ldp-mem} Let $(Z_i)_{i\in \mathbb{N}}$ be a sequence of
  $\mathbb{R}^d$-valued i.i.d. random variables where $Z_1$ satisfies
  (A-\ref{LDP-Particle}) and (A-\ref{convex-rf}). 

  Consider a triangular array $(x_i^n, 1\le i\le n,n\ge
  1)\subset {\mathcal X}$ which fulfills (A-\ref{hypo-EJP}). 

  Denote by $C_n^{\mathbf{f}}=\{ \mathbf{f}(x_i^n),\ x_i^n \in C_n\}$
  where $C_n$ is a subset of $\{x_i^n,\ 1\le i\le n\}$ given by
  Proposition \ref{proper-decomposition} and $\mathbf{f}:{\mathcal X}
  \rightarrow \mathbb{R}^{m\times d}$ is continuous. Assume that
  $C_n^{ \mathbf{f}}$ satisfies (A-\ref{Compacity}) and
  (A-\ref{Limit-Points}). Then
$$
L_n=\frac{1}{n}
\sum_1^{n} \mathbf{f}(x_i^n)\cdot Z_i
$$
satisfies the LDP in $(\mathbb{R}^m,{\mathcal B}(\mathbb{R}^m))$ with
good rate function
$$
I_{\mathbf{f}}(z) = \inf\{ \Gamma^*(z_1) + \Delta^*(z_2\mid {\mathcal
  D}),\ z_1 +z_2=z\}\ ,
$$
where the definition of ${\mathcal D}$ follows from Theorem \ref{ldp}.
\end{theo}

\begin{proof}
Recall the decomposition $L_n=\tilde L_n + \pi_n$ where 
$$
\tilde L_n =\frac 1n \sum_{x_i^n\in B_n} \mathbf{f}(x_i^n)\cdot Z_i \quad 
\textrm{and}\quad \pi_n =\frac 1n \sum_{x_i^n \in C_n} \mathbf{f}(x_i^n)\cdot Z_i,
$$
where the sets $B_n$ and $C_n$ are defined in Section
\ref{decomposition}. Theorem \ref{ldp} yields the LDP for $\pi_n$ with
good rate function $\Delta^*(\cdot \mid {\mathcal D})$. It remains now
to prove the LDP for $\tilde L_n$. We will rely on Theorem 2.2 in
\cite{Naj02} and therefore slightly modify $\tilde L_n$ so that it
fulfills the assumptions of this theorem.

In fact, it is required in \cite{Naj02} that all the points $x_i^n$ belong to ${\mathcal Y}$, 
which might not be the case here. We build in the sequel a sequence $(\tau(x_i^n))\subset{\mathcal Y}$ 
which approximates the sequence $(x_i^n,x_i^n\in B_n)$.
Let $x_i^n\in B_n$ and set
$$
\tau(x_i^n)=
\left\{
\begin{array}{ll}
x_i^n & \textrm{if}\ x_i^n\in {\mathcal Y},\\
\textrm{one of the}\ \textrm{argmin}\{ \rho(x,x_i^n),\ x\in {\mathcal Y}\} & \textrm{else}.
\end{array}
\right.
$$
Such a minimizer always exists and belongs to ${\mathcal Y}$ since ${\mathcal Y}$ 
is compact.  \\
Since $\lim_n \sup\{ \rho(x,{\mathcal Y}),\ x\in B_n\}=0$, one has $\sup_{x_i^n\in B_n} 
\rho(x_i^n,\tau(x_i^n))\xrightarrow[n\rightarrow \infty]{}0$ and 
$$
\kappa_n(\mathbf{f})\stackrel{\triangle}{=}\sup_{x_i^n\in B_n}\{
|\mathbf{f}(x_i^n)-\mathbf{f}(\tau(x_i^n))|\}\xrightarrow[n\rightarrow\infty]{}
0.
$$
Indeed, for $n$ large enough, $B_n$ lies in an $\varepsilon$-blowup of
$\mathcal Y$, which is compact since ${\mathcal X}$ is locally compact
and $\mathbf f$ is therefore uniformly continuous on this set.

Now, if we define $\bar L_n$ by
$$ \bar L_n \stackrel{\triangle}{=} \frac 1n \sum_{x_i^n\in B_n}
\mathbf{f}(\tau(x_i^n))\cdot Z_i, $$
then    $\tilde{L}_n$ and $\bar L_n$ are exponentially equivalent. Indeed,
\begin{eqnarray*}
\frac 1n \log \mathbb{P}\left( |\tilde{L}_n- \bar L_n|>\varepsilon\right)
&\le& \frac 1n \log \mathbb{P}\left( \frac 1n \sum_{i=1}^{\mathrm{card}(B_n)} |Z_i| >\frac{\varepsilon}{\kappa_n(\mathbf{f})} \right)\\
&\le& -\Lambda^*_{|Z|}\left(\frac{\varepsilon}{\kappa_n(\mathbf{f})}\right) \xrightarrow[n\rightarrow\infty]{}-\infty.
\end{eqnarray*}
where $\Lambda^*_{|Z|}$ stands for the convex conjugate of the
log-Laplace transform of $|Z|$. The measure
$\bar L_n$ satisfies all
the assumptions of Theorem 2.2 in \cite{Naj02}. Therefore, the LDP holds for it with good rate function
$\Gamma^*$. Finally the exponential equivalence yields the LDP for
$\tilde L_n$ with the same rate function 
(see for instance \cite[Theorem 4.2.13]{DemZei98}).

As the two subsums are independent, 
the contraction principle yields the LDP for 
$L_n$ with good rate function $I_{\mathbf{f}}$ given by:
\begin{equation} \label{convo}
I_{\mathbf{f}}(z) = \inf\{ \Gamma^*(z_1) + \Delta^*(z_2\mid {\mathcal D}),\ z_1 +z_2=z\}.
\end{equation}
\end{proof}

\subsection{More insight on the rate function $I_{\mathbf{f}}$}
In the convex case, that is when Assumption (A-\ref{convex-rf}) holds, the rate function $I_{\mathbf{f}}$
can be expressed more explicitely.
This section is aimed at describing how to perform the inf-convolution \eqref{convo}.

We first introduce some definitions from convex analysis (see e.g. \cite{Roc70}). 
The main result is stated in Theorem \ref{theo-rf}.

\begin{Def}[Normal cone]
Let $\mathcal C\subset \mathbb{R}^d$ be a convex set and let $a \in
\mathcal C.$ The normal cone of $\mathcal C$ at $a$, denoted by $N_{\mathcal C}(a)$, 
is defined by:
$$ 
N_{\mathcal C}(a)=\{ z\in \mathbb{R}^d;\ \langle z, x-a\rangle \le
 0,\  \forall x \in {\mathcal C}\}.
$$
\end{Def}
\begin{rem} In  particular, if $z\in N_{\mathcal C}(a)$ then $\Delta^*(z\mid
{\mathcal C}) =\langle z,a\rangle$.
\end{rem}

\begin{Def}[Relative interior]
Let ${\mathcal C}\subset \mathbb{R}^d$ be a convex set. Its affine
hull, denoted by $\mathrm{aff\,} {\mathcal C}$, is the smallest affine subset of
$\mathbb{R}^d$
containing ${\mathcal C}$. The relative interior of ${\mathcal C}$, denoted by $\mathrm{ri\,} {\mathcal C}$,
is defined by:
$$ 
\mathrm{ri\,} {\mathcal C} \stackrel{\triangle}{=} \{ x \in  \mathrm{aff\, } {\mathcal C},\ 
\exists \varepsilon >0\ \textrm{such that}\  (x + \varepsilon B(0,1)) \cap \mathrm{aff\, } \mathcal
C \subset  {\mathcal C}\}
$$
\end{Def} 
 
\begin{Def}[Subdifferential of a convex function]
A vector $x^*$ is said to be a subgradient of a convex function $f$
at a point $x$ if  for any $z$,
$$ f(z) \geq f(x) +\langle x^*, z-x \rangle.$$
The subdifferential $\partial f(x)$ of $f$ at $x$ is the set of all
subgradients
of $f$ at $x$.
\end{Def}

We can now state: 
\begin{theo}\label{theo-rf}
  Under the assumptions of Theorem \ref{ldp-mem}, the rate function $I_{\mathbf{f}}$ admits the following representation: 
\begin{equation} \label{repn1}
 I_{\mathbf{f}}(z)= \sup_{\lambda \in {\mathcal D}}
(\langle \lambda, z\rangle - \Gamma(\lambda))\ ,
\end{equation} 
where $\Gamma$ is given by \eqref{Gamma}. Furthermore, for any $ z
\in \mathrm{ri\,} \mathrm{dom}\, I_{\mathbf{f}},$ we can decompose $z$
as $z= z^* + z_{\mathbf{n}},$ where there exists $\lambda^* \in
\mathrm{dom}\, \Gamma \cap \bar{ \mathcal D}$ such that:
\begin{itemize}
\item[(i)]   $z^* \in \partial \Gamma (\lambda^*)$ and 
\item[(ii)] $z_{\mathbf{n}} \in N_{{ \bar{\mathcal D}}}(\lambda^*).$ 
\end{itemize}
In particular, for any such decomposition,
$$ I_{\mathbf{f}}(z)= \Gamma^*(z^*) + \Delta^*(z_{\mathbf{n}} \mid {\mathcal D}).$$ 
\end{theo} 

\begin{rem}[Non-exposed points]\label{non-expo} Let $z\in \mathrm{ri}\,\mathrm{dom}
  I_{\mathbf{f}}$. Consider the decomposition given by Theorem
  \ref{theo-rf}, namely $z=z^* +z_{\mathbf{n}}$, then:
  $$
  \forall t\in \mathbb{R}^+,\ I_{\mathbf{f}}(z^* + tz_{\mathbf{n}})=
  \Gamma^*(z^*) + t\langle z_{\mathbf{n}}, \lambda^*\rangle\quad
  \textrm{where}\ z^* \in \partial\Gamma(\lambda^*) \textrm{ and }
  z_{\mathbf{n}} \in N_{\mathcal D}(\lambda^*).
  $$
  In particular if $z_{\mathbf{n}} \neq 0$, $I_{\mathbf{f}}$ is affine
  in the direction $\mathbb{R}^+\ni t\mapsto z^* + tz_{\mathbf{n}}$
  and has thus infinitely many non-exposed points (see for instance
  the example developed in Section \ref{section:examples}).
\end{rem}

\begin{proof}
We first prove \eqref{repn1}. Theorem  \ref{ldp-mem} and Proposition
\ref{support-function} yield
$$ 
I_{\mathbf{f}}(z) = \inf_{z= z_1+z_2}\{\Gamma^*(z_1) +  \Delta^*(z_2\mid 
  \mathcal D)\}.
$$ 
As $ I_{\mathbf{f}}$,  $ \Gamma$ and $\Delta(.\mid \bar{
  \mathcal D})$ are convex, proper and lower semicontinuous,
we get from Theorem 16.4 in \cite{Roc70} that
\begin{eqnarray*}
  I_{\mathbf{f}}(z) &=& \left[ \Gamma +  \Delta (.\mid \bar{
      \mathcal D})\right]^*(z), \\
  &=&  \sup_{\lambda\in \mathbb R^d} \{\langle \lambda, z\rangle
  -\Gamma(\lambda) - \Delta (\lambda \mid \bar{
    \mathcal D})\}, \\
  &= & \sup_{\lambda\in \bar{\mathcal D}} \{\langle \lambda, z\rangle -\Gamma(\lambda)\} 
  \quad = \quad \sup_{\lambda\in {\mathcal D}} \{\langle \lambda, z\rangle -\Gamma(\lambda)\},
\end{eqnarray*}
and \eqref{repn1} is proved. As $ I_{\mathbf{f}}$ is convex, so is its
domain and we can consider its relative interior
$\mathrm{ri\,}\mathrm{dom\,} I_{\mathbf{f}}$. Let $z \in \ridom
I_{\mathbf{f}}$, then $I_{\mathbf{f}}(z) < +\infty$ and define $F_z$
by :
$$ 
F_z(x) = \Gamma^* (x)+ \Delta^*(z-x\mid \bar{\mathcal D}).
$$
The properties of $\Gamma^*$ and $\Delta^*(.\mid \bar{\mathcal D})$
yield that $F_z$ is proper, convex and lower semicontinuous; its level
sets are compact. In particular, the infimum of $F_z$ is attained over
$\mathbb R^d.$ Let $z^*$ be a point where this infimum is attained,
i.e.
$$
\inf_{x\in \mathbb{R}^d} F_z(x) = F_z(z^*).
$$
In this case,
$$
0 \in  \partial F_z(z^*).
$$
In order to go further in the proof, we shall describe $\partial
F_z(z^*)$ in terms of $\partial \Gamma^*$ and $\partial
\Delta^*(z-\cdot \mid \bar{\mathcal D})$. This is the purpose of the following proposition:

\begin{prop} \label{subdif}
If $z \in \mathrm{ri\, dom}\,  I_{\mathbf{f}}$ , then for any $x$,
$$ \partial F_z(x) = \partial \Gamma^*(x) - \partial \Delta^*(z-x\mid \bar{
  \mathcal D}).$$
\end{prop}   
 
\begin{proof}[Proof of Proposition \ref{subdif}]
Define $f_z$ to be the function given by $f_z(x) = \Delta^*(z-x\mid \bar{
  \mathcal D})$. Note in particular that $F_z(x)=\Gamma^*(x) +f_z(x)$.
Since $I_{\mathbf{f}}(z) = \inf_{z= z_1+z_2}\{\Gamma^*(z_1) +  \Delta^*(z_2\mid 
  \mathcal D)\},$ the sum of the epigraphs of $\Gamma^*$ and $\Delta^*$ are equal
to the epigraph of $I_{\mathbf{f}}$. This immediatly implies that 
$$ 
\mathrm{dom\,} I_{\mathbf{f}} = \mathrm{dom\,} \Gamma^* + \mathrm{dom\,} \Delta^*(\cdot \mid \bar{
  \mathcal D}).
$$
These sets being convex, Corollary 6.6.2 in \cite{Roc70} yields
$$ 
\ridom I_{\mathbf{f}} = \ridom \Gamma^* +
  \ridom \Delta^*(\cdot \mid \bar{
  \mathcal D}).$$
Let  $z \in \ridom  I_{\mathbf{f}}$, then there exists
$y \in \ridom \Gamma^*$ such that $z - y \in \ridom \Delta^*(\cdot\mid \bar{
  \mathcal D})$. This is equivalent to the fact that $y\in \ridom f_z(x)$ and therefore
\begin{equation}\label{main-assump-rocka}
\ridom \Gamma^* \cap \ridom f_z \neq \emptyset.
\end{equation}
Theorem 23.8 in \cite{Roc70} whose main assumption is fulfilled 
by \eqref{main-assump-rocka} yields then
\begin{eqnarray*}
 \partial F_z(x) & = & \partial \Gamma^*(x) + \partial f_z(x)\\
& = & \partial \Gamma^*(x)  - \partial \Delta^*(z-x\mid \bar{
  \mathcal D})
\end{eqnarray*}
and Proposition \ref{subdif} is proved.
\end{proof}
Let us now go back to the proof of Theorem
\ref{theo-rf}. By Proposition \ref{subdif},
$$
\partial F_z(z^*) = \partial \Gamma^*(z^*)  - \partial \Delta^*(z-z^* \mid \bar{
  \mathcal D}).
$$
Since $0\in \partial F_z(z^*)$, there exists $\lambda^*\in \partial \Gamma^*(z^*)$ such that 
$\lambda^* \in \partial \Delta^*(z-z^* \mid \bar{\mathcal D})$. By applying Theorem 23.5 in \cite{Roc70}, one obtains
$$
\lambda^* \in \partial \Gamma^* (z^*) \quad \Leftrightarrow \quad z^* \in \partial \Gamma(\lambda^*)
$$
which in particular implies that $\lambda^* \in \mathrm{dom\,}\Gamma$. Moreover, 
\begin{eqnarray*}
-\lambda^* \in \partial \Delta^*(z-z^*\mid \bar{\mathcal D}) &\Leftrightarrow& z-z^* \in \partial \Delta(\lambda^*\mid \bar{\mathcal D})\\
& \Leftrightarrow & z-z^* \in N_{\bar{\mathcal D}}(\lambda^*),
\end{eqnarray*} 
which in particular implies that $\lambda^* \in \bar{\mathcal D}.$\\
Denote by $z_{\mathbf{n}}= z-z^*$ , then one obtains the decomposition stated in Theorem \ref{theo-rf}. It remains to prove that:

$$ 
I_{\mathbf{f}} (z) = \Gamma^*(z^*) +  \Delta^*(z_\mathbf{n}\mid \bar{
  \mathcal D}).
$$  
We have:
\begin{eqnarray*}
I_{\mathbf{f}}(z) &=& \sup_{\lambda\in \bar{\mathcal D}} \{\langle \lambda, z\rangle -\Gamma(\lambda)\}\\
&\ge & \langle \lambda^*, z^*\rangle -\Gamma(\lambda^*) + \langle \lambda^*, z_{\mathbf{n}}\rangle 
\quad = \quad \Gamma^*(z^*) + \Delta^*(z_{\mathbf{n}} \mid {\mathcal D}) .
\end{eqnarray*}
On the other hand, 
\begin{eqnarray*}
I_{\mathbf{f}}(z) &=& \sup_{\lambda\in \bar{\mathcal D}} \{\langle \lambda, z\rangle -\Gamma(\lambda)\}\\
&\le & \sup_{\lambda\in \bar{\mathcal D}} \{\langle \lambda, z^*\rangle -\Gamma(\lambda)\}  
+\sup_{\lambda\in  \bar{\mathcal D}} \langle \lambda, z_{\mathbf{n}}\rangle 
\quad = \quad \Gamma^*(z^*) + \langle \lambda^*, z_{\mathbf{n}} \rangle,
\end{eqnarray*}
and Theorem \ref{theo-rf} is proved.
\end{proof}

\section{An example of LDP in the convex case}\label{section:examples}
To illustrate the range of Theorems \ref{ldp-mem} and \ref{theo-rf}, we
study in detail the following model :
\begin{equation} \label{def-Ln}
L_n = \frac 1 n \sum_{i=1}^n \mathbf{f}(x_i^n)  \cdot Z_i\quad \textrm{where}\quad \mathbf{f}(x)=\begin{pmatrix} 
1 & 0 \\
0 & x
\end{pmatrix}\quad \textrm{and} \quad Z_i= \begin{pmatrix} 
X_i^2 \\
X_i^2 
\end{pmatrix},
\end{equation} 
the sequence $(X_i)_{i\in \mathbb{N}}$ being a sequence of i.i.d. ${\mathcal
  N}(0,1)$ Gaussian random variables
and $(x_i^n)_{n \in \mathbb N}$ being a sequence of real numbers satisfying
$$ 
\hat R_n = \frac 1n \sum_{i=1}^n \delta_{x_i^n} \rightarrow R.
$$
We assume moreover that the support ${\mathcal Y}$ of $R$ is given by
${\mathcal Y}=[m,M]$ and that 
$$ 
\sup_{1\le i\le n} x_i^n \xrightarrow[n\rightarrow \infty]{} x_{\max} >M \qquad \textrm{and}\qquad  \inf_{1\le i\le n} x_i^n
\xrightarrow[n\rightarrow \infty]{} x_{\min} <m. 
$$
Our goal is to establish the LDP for $L_n$ and to describe as
explicitely as possible the related rate function $I_{\mathbf f}.$
\begin{rem} This example can be seen as the extension to the dimension
  2 of the example studied in \cite{BerGamRou97}. Indeed, under the same
  assumptions, Bercu et al. study the LDP for the following empirical mean
$
\frac 1 n \sum_{i=1}^n x_i^n  X^2_i\ .
$

\end{rem}
Proposition \ref{rf-ex} below is devoted to the description of the
rate function.  We first need the following notations. For $(\xi,
\xi^\prime) \in \mathbb R^2,$ set
\begin{equation}\label{Gamex}
\Gamma(\xi, \xi^\prime) = - \frac 1 2 \int \log(1-2\xi -2x
\xi^\prime) R(dx),
\end{equation}
and denote by $\Gamma^*$ the convex conjugate of $\Gamma$
(the expression for $\Gamma$ follows from a Gaussian integration and
from formula (\ref{Gamma})). Define $H$ to be the Hilbert transform of $R$, that is
$$
H(t) = \int \frac {R(dx)}{t-x} \qquad \textrm{for}\quad t\in [m,M]^c.  
$$
Set
\begin{eqnarray*}
  H_{\min} =  H(x_{\min})\qquad  &\textrm{and}&  \qquad \alpha_{\rm min} = x_{\rm min} - \frac 1{H_{\rm
      min}};\\
  H_{\rm max} =  H(x_{\max})\qquad  &\textrm{and}&  \qquad \alpha_{\rm max} = x_{\rm max} - \frac 1{H_{\rm
      max}}.
\end{eqnarray*}
Note that under the assumption that $x_{\rm min} < m$ and $x_{\rm
  max} > M$, $H_{\rm min}$ is a well-defined negative number while 
 $H_{\rm max}$ is a well-defined positive number. In particular 
$x_{\rm min} < \alpha_{\min}$ and $\alpha_{\max} < x_{\max}$. 
Moreover, the following inequalities hold true:
$$
m < \alpha_{\min} \le \int x\,R(dx)\quad \textrm{and}\quad 
\int x\, R(dx) \le \alpha_{\max} < M.
$$
In particular, $\alpha_{\min}\le \alpha_{\max}$.
In order to describe the rate function related to the LDP of $L_n$, we introduce the following domains:
\begin{eqnarray*}
  {\mathcal D}_{\infty} &=& \{(x,y)\in \mathbb{R}^2,\ x\le 0\quad \mathrm{or}\quad y \geq x_{\rm max} x\quad \mathrm{or}
\quad y  \leq x_{\min}x\} \\
  {\mathcal D}_{(I_{\mathbf{f}}=\Gamma^*)} &=& \{(x,y)\in\mathbb{R}^2,\ x> 0\quad \mathrm{and}\quad 
\alpha_{\min} x \le y\le \alpha_{\max} x \}\\
  {\mathcal D}_{\mathrm{linear}}^+ &=& \{(x,y)\in\mathbb{R}^2,\ x> 0\quad \mathrm{and}\quad 
\alpha_{\max} x < y\le x_{\max} x \}\\
  {\mathcal D}_{\mathrm{linear}}^- &=& \{(x,y)\in\mathbb{R}^2,\ x> 0\quad \mathrm{and}\quad 
x_{\min} x \le y< \alpha_{\min} x \}
\end{eqnarray*}
These domains are represented in {\sc Figure \ref{domains}} (right).
We can now state the following result.
\begin{prop}\label{rf-ex}
The empirical mean $L_n$ defined in \eqref{def-Ln} satisfies the LDP in
$\mathbb R^2$ with good rate funtion $I_{\mathbf f}$ given by
\begin{enumerate}
\item If $ (x,y)\in {\mathcal D}_{\infty}$ then $I_{\mathbf f}(x,y) = +\infty$, 
\item If $ (x,y)\in {\mathcal D}_{(I_{\mathbf{f}}=\Gamma^*)}$ then $I_{\mathbf f}(x,y) = \Gamma^*(x,y)$, 
\item If $ (x,y)\in {\mathcal D}_{\mathrm{linear}}^+$ then 
\begin{multline*}
I_{\mathbf f}(x,y) = \Gamma^*\left(H_{\rm
  max}(x_{\rm max}x-y), \alpha_{\rm max}
H_{\rm max}(x_{\rm max}x-y) \right)\\ + \frac 1 2 \left((1- H_{\rm
 max}x_{\rm max}\right)x +H_{\rm max}y),
\end{multline*}
\item If $ (x,y)\in {\mathcal D}_{\mathrm{linear}}^-$ then 
\begin{multline*}
I_{\mathbf f}(x,y) =  \Gamma^*\left(H_{\rm
 min}(x_{\rm min}x-y),  \alpha_{\rm min}H_{\rm min}(x_{\rm min}x-y) \right)\\ 
+ \frac 1 2 \left((1- H_{\rm  min}x_{\rm min})x +H_{\rm min}y\right).
\end{multline*}\\
\end{enumerate}
\end{prop}
\begin{rem}\label{non-expo-rem}
Let $x_0>0$ be fixed and consider the ray:
$$
y^-(x)=x_{\min} x +(\alpha_{\min}-x_{\min}) x_0,\quad x\ge x_0.
$$
Then 
$$
I_{\mathbf{f}}(x,y^-(x))=\Gamma^*(x_0,\alpha_{\min} x_0) +\frac 12 (x-x_0).
$$
In particular, there are infinitely many non-exposed points for
$I_{\mathbf{f}}$ along the ray $((x,y^-(x)); x\ge x_0)$.  The same can
be shown along the ray
$$
y^+(x)=x_{\max} x +(\alpha_{\max} -x_{\max}) x_0;\ x\ge x_0.
$$
\end{rem}

\begin{proof}[Proof of Proposition \ref{rf-ex}]
The LDP will be established as soon as assumptions of Theorem \ref{ldp-mem} are fulfilled. 
It is straightforward to check (A-\ref{LDP-Particle}) to
(A-\ref{Compacity}) and (A-\ref{hypo-EJP}). In order to check Assumption (A-\ref{Limit-Points}),
we rely on the following lemma: 
\begin{lemma}\label{Dmax} For every $x\in [x_{\rm min}, x_{\rm max}]$, one has:
$$ 
\mathcal D_{\mathbf f (x_{\rm min})} \cap  \mathcal D_{\mathbf f
  (x_{\rm max})} \subset \mathcal D_{\mathbf f (x)}.$$ 
\end{lemma}
\begin{proof}[Proof of Lemma \ref{Dmax}]
Let $(\xi, \xi^\prime) \in \mathcal D_{\mathbf f (x_{\rm min})} \cap
\mathcal D_{\mathbf f (x_{\rm max})}.$
This implies that  $(\xi, x_{\rm min}\xi^\prime) \in \mathcal D_{Z_1}$ and 
$(\xi, x_{\rm max}\xi^\prime) \in \mathcal D_{Z_1}.$
Every $x \in [x_{\rm min}, x_{\rm max}]$ can be written as a convex
combination
of $x_{\rm min}$ and $ x_{\rm max}:$
$ x= a x_{\rm min} + b x_{\rm max},$ where $a+b=1$, $a,b$ being nonnegative.
By convexity of  $\mathcal D_{Z_1}$,
$ (\xi, x \xi^\prime) =  
 a (\xi, x_{\rm min}\xi^\prime) + b (\xi, x_{\rm max}\xi^\prime)\in
 \mathcal D_{Z_1}.$
Therefore $(\xi, \xi^\prime) \in \mathcal D_{\mathbf f(x)}.$ 
\end{proof}
We can now check (A-\ref{Limit-Points}). The mere definition of $x_{\min}$ and $x_{\max}$ 
implies that both $x_{\min}$ and $x_{\max}$ belong to $\Cout^{\mathbf{f}}$ and $\Cin^{\mathbf{f}}$ and that
both $\Cout^{\mathbf{f}}$ and $\Cin^{\mathbf{f}}$ are included in $[x_{\min},x_{\max}]$.
In particular, the set $\mathcal D$ is well defined and is given by:
$$
{\mathcal D}= \bigcap_{\{x,\ \mathbf{f}(x) \in
  \Cout^{\mathbf{f}}\}}{\mathcal D}_{\mathbf{f}(x)}
\stackrel{(a)}{=}\mathcal D_{\mathbf f (x_{\rm min})} \cap \mathcal
D_{\mathbf f (x_{\rm max})}\stackrel{(b)}{=} \bigcap_{\{x,\
  \mathbf{f}(x) \in \Cin^{\mathbf{f}}\}} {\mathcal D}_{\mathbf{f}(x)}
$$
where $(a)$ and $(b)$  follow from Lemma \ref{Dmax}. An easy computation yields 
\begin{equation} \label{D}
 \mathcal D = \{(\xi, \xi^\prime) \in \mathbb R^2;\ 
1-2\xi-2x_{\rm min}\xi^\prime >0 \textrm{ and }  1-2\xi-2x_{\rm max}\xi^\prime >0\}.
\end{equation}
The LDP is therefore established by applying Theorem 
\ref{ldp-mem} and the rate function is given by:
$$ 
I_{\mathbf f }(z) = \inf_{z=z_1+z_2} \{\Gamma^*(z_1) +
\Delta^*(z_2| \mathcal D)\},
$$
with  $\mathcal D$ as above and $\Gamma$ as defined in \eqref{Gamma}.
Formula \eqref{Gamex} yields: 
$$ 
\textrm{dom}\, \Gamma = \{(\xi, \xi^\prime) \in \mathbb R^2;\ 
1-2\xi-2x\xi^\prime >0 \textrm{ for all }  x \in [m,M]\},
$$
and therefore
\begin{equation} \label{domG}
\textrm{dom}\, \Gamma = \{(\xi, \xi^\prime) \in \mathbb R^2;\ 
1-2\xi-2m\xi^\prime >0 \textrm{ and }  1-2\xi-2M\xi^\prime >0\}.
\end{equation}
{\sc Figure} \ref{figex} shows $\textrm{dom}\, \Gamma$ and $\mathcal D $
for particular choices of the parameters.\\

\begin{figure}
\begin{center}
\epsfig{figure=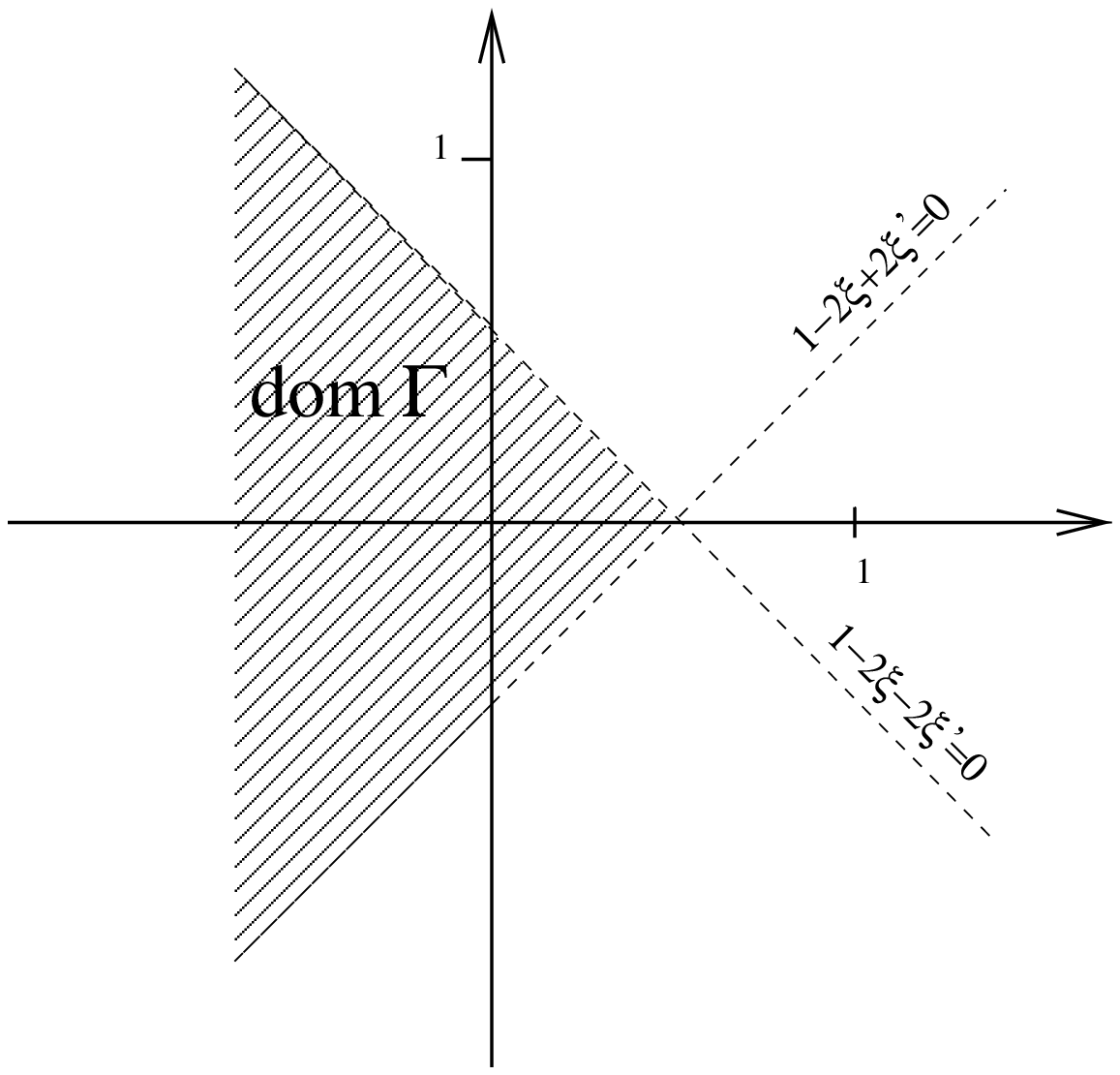,height=5cm}
\epsfig{figure=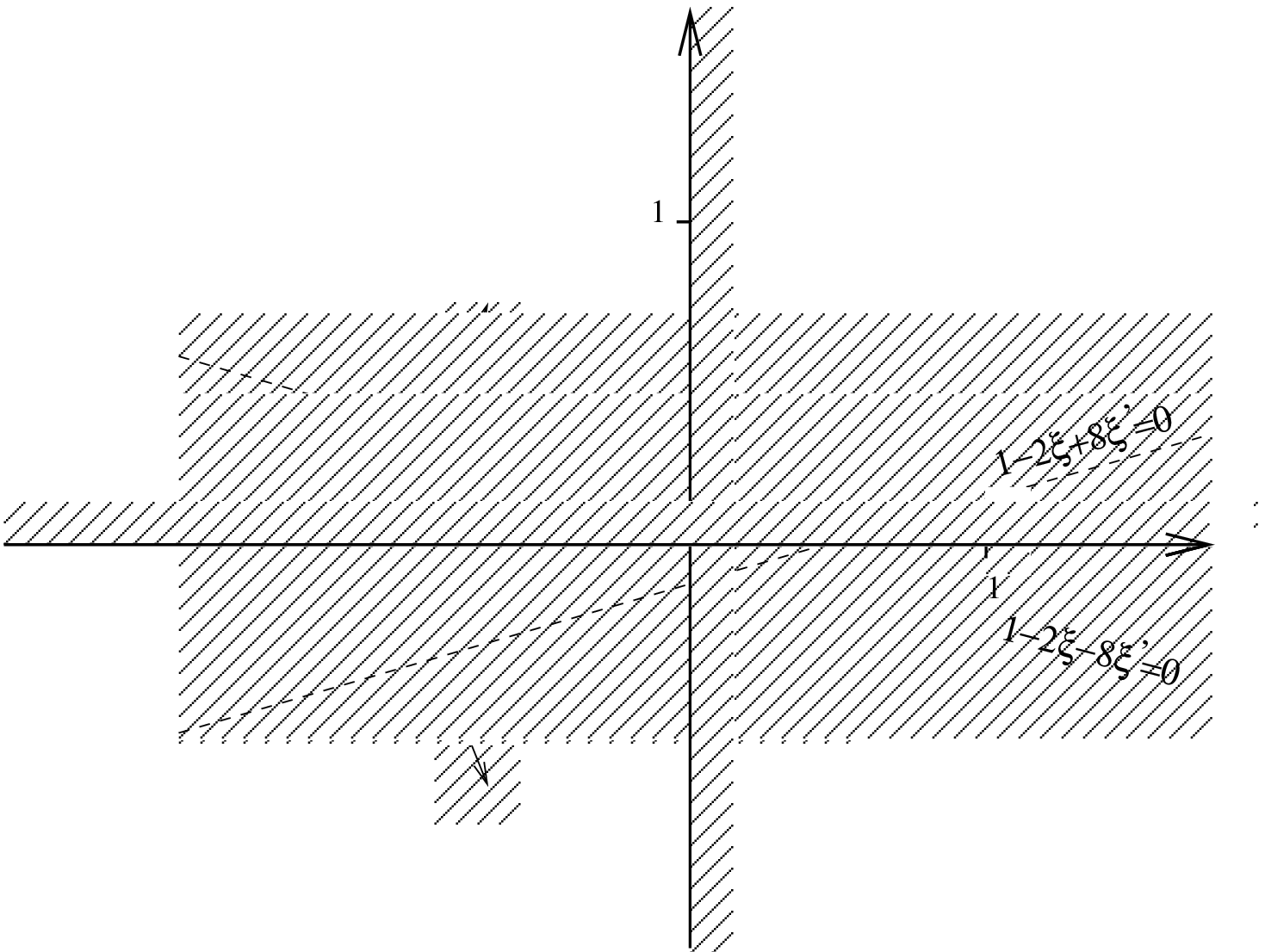,height=5cm}
\caption{\small On this figure are represented $\textrm{dom}\,
  \Gamma$ 
for $m=-1$ and $M=1$ (left) and
  $\mathcal D $
for $x_{\rm min} =-4$ and $x_{\rm max} =4$ (right).
On the picture of $\mathcal D,$ we figured also 
some of the normal cones to $\bar{\mathcal D}$, whose directions are represented by the arrows.}\label{figex}
\end{center}
\end{figure} 

We first prove Proposition \ref{rf-ex}-(1). In order to prove this statement, it is equivalent 
to determine the domain of $I_{\mathbf f}.$ We use the fact that 
$$ 
\textrm{dom}\,I_{\mathbf f} = \textrm{dom}\, \Gamma^* + 
\textrm{dom}\, \Delta^*(\cdot \mid \mathcal D) 
$$
and focus on the two domains of the right-hand side. One can check that
\begin{eqnarray*} 
\textrm{dom}\, \Gamma^* &=& \{(x, y) \in \mathbb R^2;\ x>0 \textrm{ and } 
mx\leq y \leq Mx\},\\
\textrm{dom}\, \Delta^*(\cdot\mid\mathcal D) &=& \{(x, y) \in \mathbb R^2;\ 
x\geq 0 \textrm{ and } x_{\rm min} x \leq y\leq x_{\rm max}
x\}.
\end{eqnarray*}
Therefore
 \begin{equation} 
\textrm{ dom}\, I_{\mathbf f} = \{(x, y) \in \mathbb R^2;\ 
x>0\quad \textrm{and}\quad 
x_{\rm min} x < y < x_{\rm max}x\}.
\end{equation}
Note in particular that in this case, $\textrm{ ri
dom}I_{\mathbf f} =\textrm{ 
dom}I_{\mathbf f} . $\\

The three domains $\textrm{dom}\, \Gamma^*,$ $\textrm{dom}\, \Delta^*(\cdot\mid\mathcal D)$
and $\textrm{dom}\,  I_{\mathbf f} $ are represented on {\sc Figure} \ref{domains}.\\ 

\begin{figure}
\begin{center}
\epsfig{figure=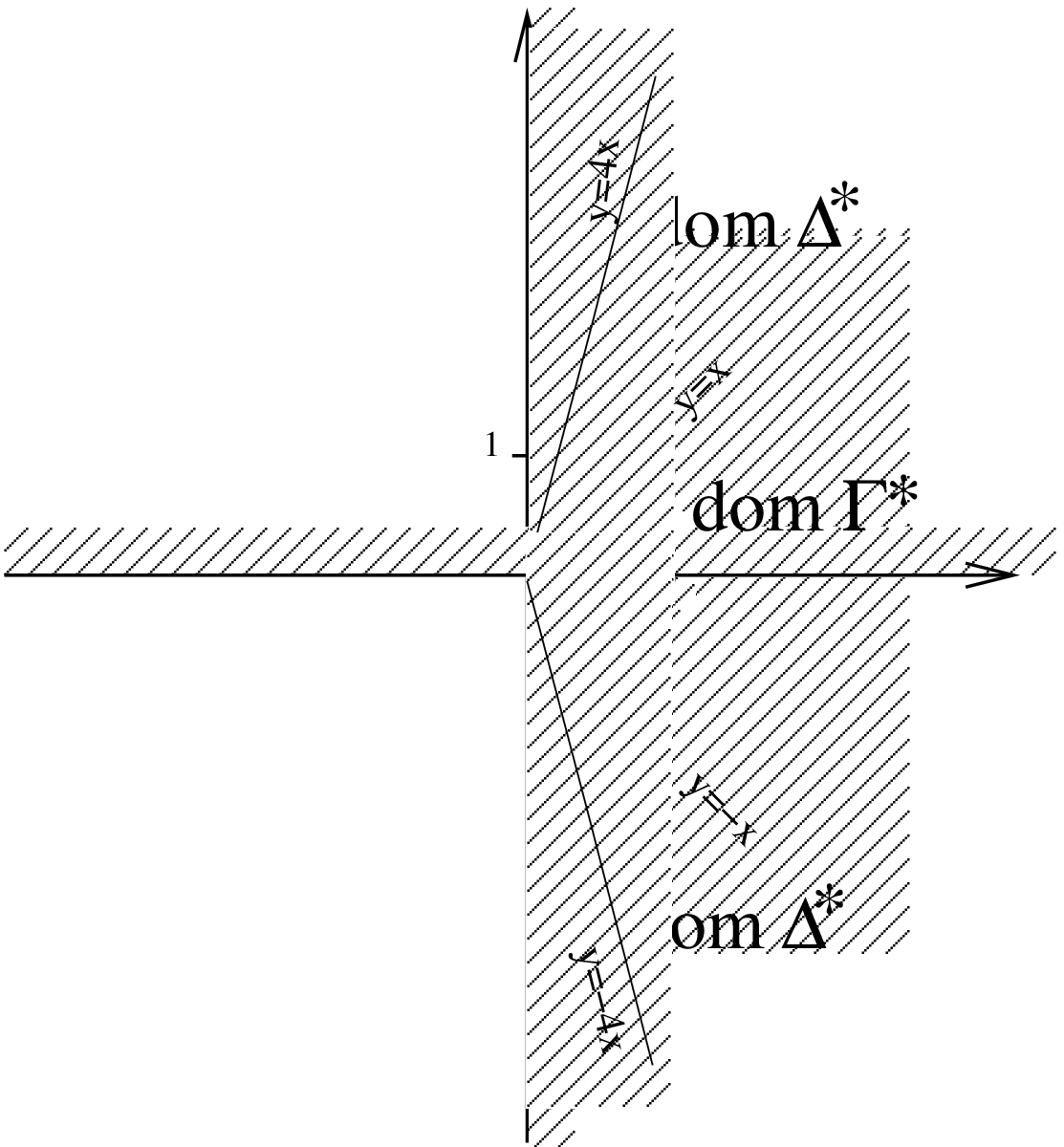,height=5cm}
\epsfig{figure=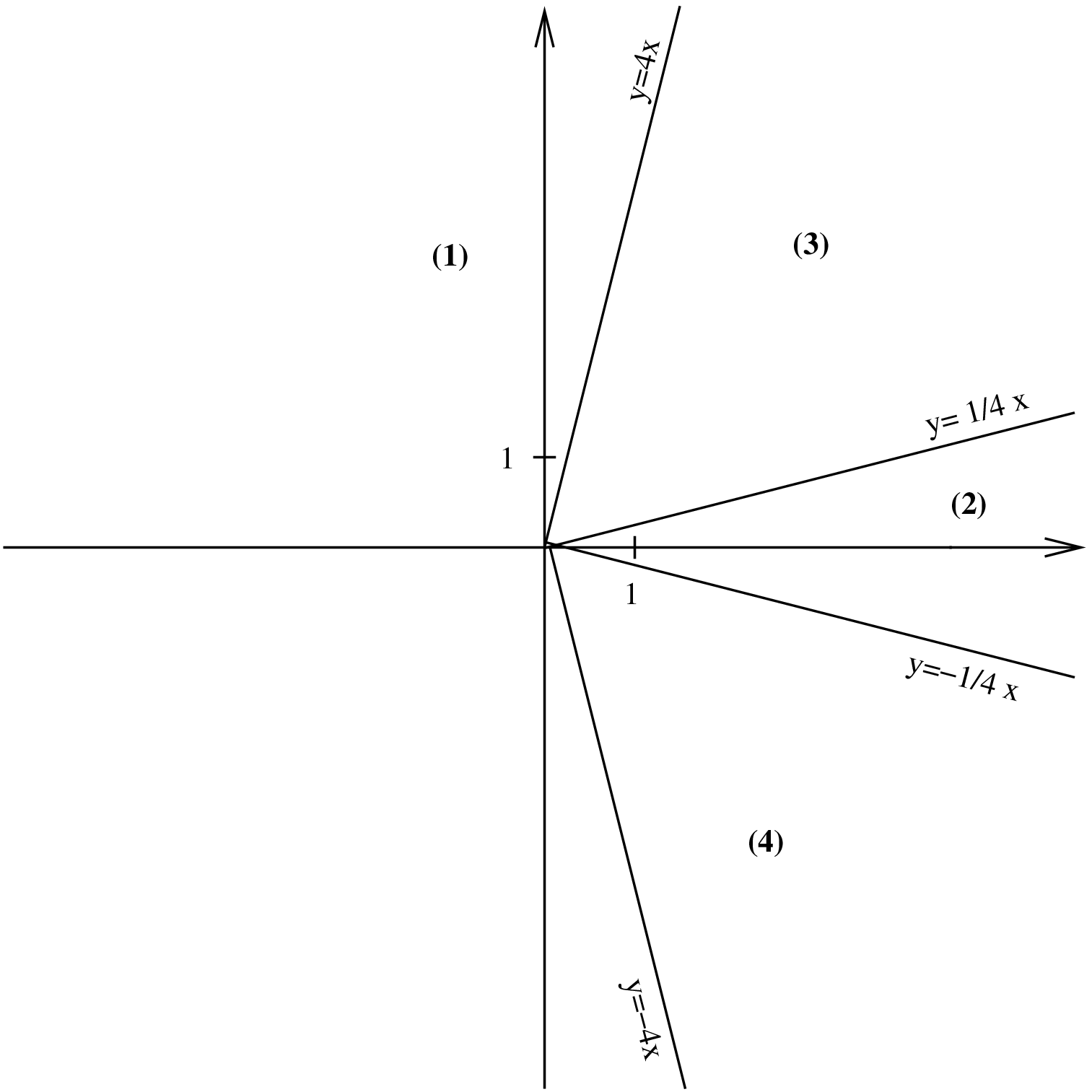,height=5cm}
\caption{\small The left picture represents 
$\textrm{dom}\, \Gamma^*$ (hatched cone) and  $\textrm{dom}\, \Delta^*(\cdot\mid
\mathcal D)$
(delimited by the two half-lines $y=4x$ and $y=-4x$).
The right picture represents the four zones of $\mathbb{R}^2$ where 
$I_{\mathbf{f}}$ has a particular expression. Zone (1) (resp. (2), (3) and (4)) 
represents ${\mathcal D}_{\infty}$ (resp. ${\mathcal D}_{(I_{\mathbf{f}}=\Gamma^*)}$,
${\mathcal D}_{\mathrm{linear}}^+$ and ${\mathcal D}_{\mathrm{linear}}^-$).
We kept the same values
of the parameters as in {\sc Figure} \ref{figex} and chose a
particular $R$ for which
$H_{\rm max} = - H_{\rm min} = 4/15$.}\label{domains}
\end{center}
\end{figure} 

We now prove Proposition \ref{rf-ex}-(2).
Theorem \ref{theo-rf} yields:
$$ 
I_{\mathbf f}  (z) = \sup_{\lambda \in \bar{\mathcal D}}\{\la
\lambda, z \ra - \Gamma(\lambda)\}.
$$
If one consider $g_z(\lambda) = \la \lambda, z \ra - \Gamma(\lambda),$
one can check that for $z \in \textrm{dom}\, \Gamma^*,$ an element
$\bar \lambda = (\bar \xi, \bar \xi^\prime)$  realizing the supremum
of $g_z$ satisfies the condition
$$ 
\alpha - \frac 1{H(\alpha)} = \frac y x, \quad \textrm{with } \alpha = \frac{1-
  2\bar \xi}{2  \bar \xi^\prime}. 
$$
Therefore $\bar \lambda \in  \textrm{dom}\, \Gamma \cap \bar{\mathcal
  D}$  if and only if  $\frac y x \in [\alpha_{\rm min}, \alpha_{\rm max}]$
and in this case $I_{\mathbf f} (z)= \Gamma^*(z). $\\

We now turn to the proof of Proposition \ref{rf-ex}-(3).
>From Theorem \ref{theo-rf}, we just need to exhibit a decomposition
$z= z^* + z_{\mathbf{n}},$ where $z^* \in \partial\Gamma(\lambda^*)$
and $z_{\mathbf{n}} \in N_{\bar{\mathcal D}}(\lambda^*)$ for some $
\lambda^* \in  \textrm{dom} \Gamma \cap \bar{\mathcal D}$. 
In this case, the value of $I_{\mathbf f} (z)$ is given by $I_{\mathbf f} (z) = \Gamma^*(z^*) 
+\langle \lambda^*, z_{\mathbf{n}}\rangle$.
One can check that $ \textrm{dom}\, \Gamma \cap \bar{\mathcal
  D}$ can be split into three subsets : the interior of $\mathcal D$,
and the two half-lines $\{1-2\xi-2x_{\rm min} \xi^\prime=0, \xi < 1/2\}$
and $\{1-2\xi-2x_{\rm max} \xi^\prime=0, \xi < 1/2\}$. The normal cones
to $ \bar{\mathcal D}$ are then easy to determine:
\begin{itemize}
\item[-] if $(\xi, \xi^\prime) \in \textrm{int } \mathcal D$, then
$N_{\bar{\mathcal D}}(\xi, \xi^\prime) = \{(0,0)\},$
\item[-] if $\xi< 1/2$ 
and $ 1-2\xi-2x_{\rm min}\xi^\prime =0 ,$ then
$N_{\bar{\mathcal D}}(\xi, \xi^\prime) = \{t(1,x_{\rm min}), t\geq0\},$
 \item[-] if $\xi < 1/2$ 
and $ 1-2\xi-2x_{\rm max}\xi^\prime =0 ,$ then
$N_{\bar{\mathcal D}}(\xi, \xi^\prime) = \{t(1,x_{\rm max}), t\geq0\}.$
\end{itemize}
These normal cones are represented by the arrows on {\sc Figure \ref{figex}}(right).

\begin{figure}[htbp]\label{decompo}
\begin{center}
\epsfig{figure=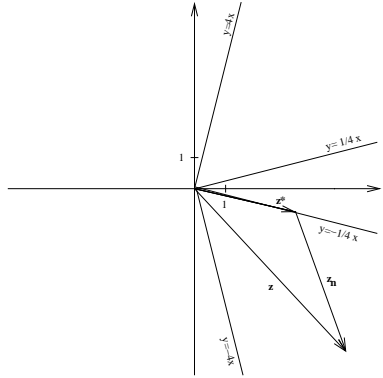,height=5cm}
\caption{\small For a $z =(x,y)$ such that  $x_{\rm min} x< y< \alpha_{\rm
    min} x,$ 
we decompose  $z = z^* + z_n$ with $z^*$ such that 
$y^* = \alpha_{\rm min} x^*$  and $z_n = t(1, x_{\rm min} ),$
for a $t>0$.}
\end{center}
\end{figure}

We can now conclude the proof of the third point of the proposition.
If we choose 
\begin{eqnarray*}
\lambda^* &=& \left(\frac 1 2 - \frac{x_{\rm
      min}}{y-x_{\rm min }x}, \frac 1{y-x_{\rm min}x}\right), \\
z^* &=&
(H_{\rm min}(x_{\rm min}x-y), (x_{\rm min}H_{\rm min}-1)(x_{\rm
  min}x-y)),\\
z_{\mathbf{n}}&=&z-z^*,
\end{eqnarray*}
it is easy to check that
this decomposition fulfills the required properties, i.e.  $z^* \in
\partial\Gamma(\lambda^*)$ and $z_{\mathbf{n}} \in N_{\bar{\mathcal
    D}}(\lambda^*)$ for some $ \lambda^* \in \textrm{dom} \Gamma \cap
\bar{\mathcal D}$. Therefore,
\begin{eqnarray*}
 I_{\mathbf f}(z) & = & \Gamma^*(z^*) + \la \lambda^*, z_n \ra \\
 & = & \Gamma^*(z^*) + \frac 1 2 (x + H_{\rm min}(y-x_{\rm min}x))
\end{eqnarray*}
The decomposition $z= z^* + z_{\mathbf{n}}$ 
can be seen on Figure \ref{decompo}.

The proof of Proposition \ref{rf-ex}-(4) is very similar and is left to the reader.
\end{proof}

\subsection*{Remarks on the LDP and the spherical integral}
We conclude this section with remarks related to the prime motivation
of this study, namely the study of the asymptotics of spherical
integrals. We recall from \cite{GuiMai05} that the goal is to get the
asymptotics of
\begin{equation}\label{spherical-integral}
I_n(A_n,B_n) = \int e^{N\ {\rm Trace} (A_nUB_nU^*)} dm_n(U),
\end{equation}
where $A_n$ and $B_n$ are two real diagonal matrices and $m_n$ is the
Haar measure on the orthogonal group.  Obtaining the asymptotic
expansion of such integrals has major applications in statistics for
instance.  Indeed, the asymptotic expansion for the joint eigenvalue
density of some deformed Wigner matrices can readily be deduced
from the above integral.

In the case where $A_n$ is of rank one, with a unique nonzero eigenvalue
denoted by $\theta$ and where $B_n=\mathrm{diag}(x_i^n,\ 1\le i\le n)$ where 
$\frac 1n \sum \delta_{x_i^n}$ converges, the spherical integral can be written as 
\begin{equation}
\label{IN}
I_n(A_n,B_n) = \mathbb E\,\exp\left(n\theta \frac{\sum_{i=1}^n x_i^n
    X_i^2}{\sum_{i=1}^n X_i^2}\right),
\end{equation}
where $\mathbb E$ is the expectation under the standard
$N$-dimensional Gaussian measure.

A natural strategy to tackle the asymptotics of $I_n$ is then to establish the LDP for the empirical 
measure $L_n$ as studied in the previous example and to apply Varadhan's lemma
to get the asymptotics of $I_n$ (see \cite[Theorem 6]{GuiMai05}).

Beside the fact that we fully recover the LDP result of
\cite{GuiMai05}, we believe that the representation of the rate
function (Theorem \ref{theo-rf}) sheds new light on the role played by
the largest and lowest eigenvalues in the asymptotics of the rank-one
spherical integral: The very reason comes from the fact that the individual rate
function of the particle $\frac 1n {\tiny \left( \begin{array}{c} X_1^2\\ X_1^2 \end{array}\right)}$
fulfills the convexity assumption (A-\ref{convex-rf}). This is in particular illustrated 
in Lemma \ref{Dmax}.

In the forthcoming section, we study the LDP in the non-convex case,
that is when (A-\ref{convex-rf}) is not fulfilled. This will lead to
partial results in the study of the asymptotics of the spherical
integral beyond the rank-one case.


\section{The LDP in the non-convex case}\label{section:nonconvex}

There are several models which fulfill Assumption
(A-\ref{LDP-Particle}) with a non-convex rate function.  Take for
instance the simple model $Z_1=(X_1^2, Y_1^2, X_1 Y_1)$ where $X_1$
and $Y_1$ are independent standard Gaussian random variables. Denote
by ${\mathcal C}=\{(x,y,z) \in \mathbb R^3,\ z=-\sqrt{xy} \textrm{ or }
z= \sqrt{xy}\}$, then $\frac{Z_1}n $
satisfies the LDP with good rate function
$$
I(x,y,z)=\frac x2 +\frac y2 +\Delta (z \mid {\mathcal C})\quad \textrm{where}\quad 
\Delta(z\mid {\mathcal C})=\left\{
\begin{array}{ll}
0& \textrm{if}\ z\in {\mathcal C}\\
\infty & \textrm{else}.
\end{array}
\right.,
$$
which is highly non-convex. We will see that this kind of models arises in the study of
spherical integrals and may give rise to interesting phenomenas. 

We give in this section an assumption over the set $A_n=\{x_i^n\in
{\mathcal X},\,1\le i\le n\}$ which ensures the LDP for $L_n$ to hold.
Although quite stringent, this assumption encompasses interesting
models as we shall see. We then state the LDP.

Recall that ${\mathcal Y}$ is the support of the limiting probability $R$.
\begin{assump}\label{assumption-nonconvex}
Assume that ${\mathcal X}\subset \mathbb{R}^p$ for a given integer $p$. 
Denote by $A_n=\{x_i^n\in {\mathcal X},\ 1\le i\le n\}$. 
Then there exists an integer $T$ such that:
$$
A_n=\tilde A_n \cup \bigcup_{\ell=1}^T \{x_{i_{\ell}}^n\}
$$
where $\rho(\tilde A_n,{\mathcal Y})$ goes to zero as $n\rightarrow \infty$ 
while for $1 \leq \ell\leq T$,
$$
x_{i_{\ell}}^n \xrightarrow[n\rightarrow \infty]{} x_{\ell}^{\infty} ,
$$
where the $x_{\ell}^{\infty}$'s do not belong to ${\mathcal Y}$.
\end{assump}

\begin{rem} Assumption (A-\ref{assumption-nonconvex}) implies that there 
exists a finite number of outliers $x_{i_{\ell}}^n$ that remain outside 
the support ${\mathcal Y}$ and that converge pointwise to a limit
$x_{\ell}^\infty$.
\end{rem}

\begin{theo}\label{ldp-nonconvex}
Assume that $(Z_i)_{i\in \mathbb{N}}$ is a sequence of
  $\mathbb{R}^d$-valued i.i.d random variables where $Z_1$  satisfies
  (A-\ref{LDP-Particle}). Assume that (A-\ref{hypo-EJP}) and (A-\ref{assumption-nonconvex}) hold 
for the sequence $(x_i^n, 1\le i\le n,n\ge 1)$. Then 
$$
L_n=\frac{1}{n}
\sum_1^{n} \mathbf{f}(x_i^n)\cdot Z_i
$$
satisfies the LDP in $(\mathbb{R}^m,{\mathcal B}(\mathbb{R}^m))$ with
good rate function
\begin{equation}\label{rf-nonconvex}
I_{\mathbf{f}}(z) = \inf\left\{ \Gamma^*(z_0) + \sum_{\ell=1}^T I(y_{\ell}) ;\ 
z_0+\sum_{\ell=1}^T \mathbf{f}(x_{\ell}^\infty)\cdot y_{\ell} = z\right\}.
\end{equation}
\end{theo}

\begin{proof}
Recall that $A_n=\tilde A_n \cup \bigcup_{\ell =1}^T \{x_{i_{\ell}}^n\}$ 
by (A-\ref{assumption-nonconvex}) and write:
$$
L_n= \frac 1n \sum_{x_i^n\in \tilde A_n} \mathbf{f}(x_i^n)\cdot Z_i + \frac 1n 
\sum_{\ell=1}^T \mathbf{f}(x_{i_{\ell}}^n)\cdot Z_{i_{\ell}},
$$
One can prove the LDP for $\frac 1n \sum_{x_i^n\in \tilde A_n}
\mathbf{f}(x_i^n)\cdot Z_i$ as in the proof of Theorem \ref{ldp-mem}
(which relies on an adaptation of Theorem 2.1 in \cite{Naj02} and does
not involve the convexity of $I$). On the other hand,
$\sum_{\ell=1}^T \frac{\mathbf{f}(x_{i_{\ell}}^n)\cdot Z_{i_{\ell}}}{n}$ is
exponentially equivalent to $\sum_{\ell=1}^T \frac{\mathbf{f}(x_{\ell}^\infty)\cdot
  Z_{i_{\ell}}}{n}$ which satisfies the LDP with good rate function 
$$
J(z)=\inf\left \{ \sum_{\ell=1}^T I(y_{\ell}),\ \sum_{\ell=1}^T \mathbf{f}(x_{\ell}^\infty)\cdot y_{\ell}=z\right\}.
$$
Since $\frac 1n \sum_{x_i^n\in \tilde A_n} \mathbf{f}(x_i^n)\cdot Z_i$ and $\frac 1n 
\sum_{\ell=1}^T \mathbf{f}(x_{i_{\ell}}^n)\cdot Z_{i_{\ell}}$ are independent, the LDP holds
with good rate function $I_{\mathbf{f}}$ given by \eqref{rf-nonconvex}. Proof of Theorem \ref{ldp-nonconvex}
is completed.
\end{proof}

\section{An example of LDP in the non-convex case: Influence of the
  second largest eigenvalue}\label{section:eigen}

\subsection{Presentation of the example}
In this section, we shall study a simple model which underlines the
differences between the LDP in the convex case and the LDP in the non-convex one.
Consider the set $A_n=\{x_i^n,\ 1\le i\le n\}$ where $x_1^n=\kappa_1$, $x_2^n=\kappa_2$ and $x_i^n=1$ for $i\ge 3$.
Assume the following:
$$
1<\kappa_2<\kappa_1.
$$
One can think of the $x_i^n$ as the eigenvalues of a $n\times n$
matrix and one can check that 
$$
\frac 1n \sum_{i=1}^n \delta_{x_i^n} \xrightarrow[n\rightarrow\infty]{} \delta_1
$$ 
while $\kappa_1$ and $\kappa_2$ are two outliers.

In the sequel, we study the influence of the second
largest eigenvalue $\kappa_2$ over the rate function of a given LDP in a
convex and non-convex case. We prove that the second largest eigenvalue has no influence 
on the rate function that drives the LDP in the convex case (Proposition \ref{second-eigen-convex}) 
while this eigenvalue has an impact 
on the LDP in the non-convex case (Proposition \ref{second-eigen-nc}). We finally go back 
to spherical integrals and make some concluding remarks.

Denote by $\mathbf{f}$ the following matrix-valued function:
$$
\mathbf{f}(x)=\left({\tiny \begin{array}{ccc}
1 & 0 & 0\\
0 &1 &0 \\
x &0 &0 \\
0 &x &0 \\
0 &0 &1
\end{array}}\right)
$$ 

Let us now introduce the random variables we will consider.

\subsection{The convex model} Consider a family of $\mathbb{R}^3$-valued random variables 
$( Z_i)_{i\ge 1}$ satisfying Assumptions (A-\ref{LDP-Particle}) and (A-\ref{convex-rf}).
Denote by
\begin{eqnarray*}
L_n( Z) &=& \frac 1n \sum_{i=1}^n \mathbf{f}(x_i^n)\cdot  Z_i\\
&=& \frac 1n \mathbf{f}(\kappa_1)\cdot  Z_1 
+\frac 1n \mathbf{f}(\kappa_2)\cdot  Z_2 
+\frac 1n \sum_{i=3}^n \mathbf{f}(x_i^n)\cdot  Z_i\\
&\stackrel{\triangle}{=}& \pi_n^1( Z) + \pi_n^2( Z) +\tilde{L}_n( Z)\\
\textrm{and by}\ \bar L_n( Z)&\stackrel{\triangle}{=}& \pi_n^1( Z)+\tilde{L}_n( Z)\\
\end{eqnarray*}

One can apply Theorem \ref{ldp-mem} to $L_n( Z)$ and $\bar
L_n( Z)$ which therefore satisfy LDPs with given rate functions that
we denote respectively by $I_{
  Z}$ and $\bar I_{ Z}$.

\begin{prop}\label{second-eigen-convex}
The rate functions $I_{ Z}$ and $\bar I_{ Z}$ related to the LDPs of $L_n( Z)$ and 
$\bar L_n( Z)$ are equal.
\end{prop}
\begin{rem} This proposition underlines the fact that the second largest eigenvalue 
does not have any influence on the rate function of the LDP.
\end{rem}

\begin{proof} Let 
$$
 Z_i= \left(
\begin{array}{c}
U_i\\
V_i\\
W_i
\end{array}\right)
\qquad \textrm{then}\qquad 
\mathbf{f}(x)\cdot  Z_i =\left({\tiny 
\begin{array}{c}
U_i\\
V_i\\
xU_i\\
xV_i\\
W_i
\end{array}}\right).
$$
For $\lambda \in \mathbb R^5$, denote by
\begin{eqnarray*}
\Lambda(\lambda) &=& \ln \mathbb{E} e^{\langle \lambda,\mathbf{f}(1) \cdot  Z \rangle},\\
\Lambda_i(\lambda) &=& \ln \mathbb{E} e^{\langle \lambda, \mathbf{f}(\kappa_i)\cdot  Z\rangle },\quad i\in \{1,2\}.
\end{eqnarray*}
Consider also the associated domains:
\begin{eqnarray*}
{\mathcal D}_0&=&\{ \lambda \in \mathbb{R}^5;\ \Lambda(\lambda)<\infty\},\\
{\mathcal D}_i&=&\{ \lambda \in \mathbb{R}^5;\ \Lambda_i(\lambda)<\infty\},\quad i\in \{1,2\}.
\end{eqnarray*}
Remark that 
\begin{equation}\label{convex-domain}
\lambda = (\alpha,\beta,\gamma,\delta,\theta)\in {\mathcal D}_i \quad \Leftrightarrow\quad 
\lambda_i = (\alpha,\beta,\kappa_i \gamma,\kappa_i \delta,\theta)\in {\mathcal D}_0, \quad i\in \{1,2\}.
\end{equation}
>From Theorem \ref{ldp-mem}, we know that 
$$
  I_{ Z}(z) = \sup_{\lambda\in {\mathcal D}_0\cap {\mathcal D}_1 \cap {\mathcal D}_2} 
\{\langle \lambda,z\rangle -\Lambda(\lambda)\}\quad \textrm{and}\quad 
 \bar{I}_{ Z}(z) = \sup_{\lambda\in {\mathcal D}_0\cap {\mathcal D}_1}
\{\langle \lambda,z\rangle -\Lambda(\lambda)\}
$$
We now prove that 
$
\lambda\in {\mathcal D}_0\cap {\mathcal D}_1$ implies that $\lambda\in {\mathcal D}_2.$ Let $\lambda =(\alpha,\beta,\gamma,\delta,\theta)\in{\mathcal D}_0\cap {\mathcal D}_1$.
>From \eqref{convex-domain},
$$
\lambda \in {\mathcal D}_1 \quad \Rightarrow \quad
\lambda_1=(\alpha,\beta,\kappa_1 \gamma,\kappa_1 \delta,\theta)\in
{\mathcal D}_0.
$$   
Moreover, as $1< \kappa_2 < \kappa_1$, $\kappa_2 $ can be written as 
$\kappa_2  = a + b \kappa_1,$ with  $a,b$ non-negative and
$a+b=1.$
Due to the convexity of ${\mathcal D}_0$, we have that
$
a \lambda +b \lambda_1  \in {\mathcal D}_0.
$
On the other hand, 
$$ a \lambda +b \lambda_1 = (\alpha, \beta, \kappa_2 \gamma, \kappa_2 \delta,
\theta),$$
so that  $\lambda\in {\mathcal D}_2$ by \eqref{convex-domain}.
Therefore,
\begin{eqnarray*}
I_{ Z}(z) &=& \sup_{\lambda\in {\mathcal D}_0\cap {\mathcal D}_1 \cap {\mathcal D}_2} 
\{\langle \lambda,z\rangle -\Lambda(\lambda)\}\\
&=&\sup_{\lambda\in {\mathcal D}_0\cap {\mathcal D}_1}
\{\langle \lambda,z\rangle -\Lambda(\lambda)\} \quad =\quad \bar{I}_{ Z}(z) 
\end{eqnarray*}
and the proof of Proposition \ref{second-eigen-convex} is completed.

\end{proof}

\subsection{The non-convex model}\label{example-non-convex-sub}
Let $(X_i)_{i\ge 1}$ and $(Y_i)_{i\ge 1}$ be two independent families
of i.i.d. standard Gaussian random variables and consider the i.i.d. 
$\mathbb{R}^3$-valued 
random variables 
\begin{eqnarray*}
{\check Z}_i=\left(
\begin{array}{c}
X_i^2\\
Y_i^2\\
X_i Y_i
\end{array}
\right).
\end{eqnarray*}
We shall study the LDP of 
\begin{eqnarray*}
L_n({\check Z})&=&\frac 1n \sum_{i=1}^n \mathbf{f}(x_i^n)\cdot {\check Z}_i\\
&=&\frac 1n \left( {\tiny \begin{array}{c}
X_1^2\\
Y_1^2\\
\kappa_1 X_1^2\\
\kappa_1 Y_1^2\\
X_1 Y_1 
\end{array}}\right) +
  \frac 1n 
\left( {\tiny \begin{array}{c}
X_2^2\\
Y_2^2\\
\kappa_2 X_2^2\\
\kappa_2 Y_2^2\\
X_2 Y_2 
\end{array}}\right) +
\frac 1n \sum_{i=3}^n
\left( {\tiny \begin{array}{c}
X_i^2\\
Y_i^2\\
X_i^2\\
Y_i^2\\
X_i Y_i 
\end{array}}\right) \\
&\stackrel{\triangle}{=}& \pi^1_n({\check Z}) + \pi_n^2({\check Z}) +\tilde{L}_n({\check Z})
\end{eqnarray*}
As above, we also introduce $\bar{L}_n({\check Z})=\pi_n^1({\check Z})+\tilde{L}_n({\check Z})$.\\
\\
\indent The non-convex model satisfies assumptions of Theorem \ref{ldp-nonconvex}. Therefore, both 
$L_n({\check Z})$ and $\bar{L}_n({\check Z})$ satisfy the LDP with
given rate functions 
that we denote respectively by $I_{\check Z}$ and $\bar I_{\check Z}$. 

We shall prove the following:
\begin{prop}\label{second-eigen-nc}
  Let $\kappa_1 <2\kappa_2 -1$. The rate function $I_{\check Z}$ that drives the LDP for $L_n({\check Z})$ differs
  from the rate function $\bar I_{\check Z}$ that
  drives the LDP for $\bar{L}_n({\check Z})$.
\end{prop}

\begin{rem}
  Proposition \ref{second-eigen-nc} illustrates the influence
  of the second largest eigenvalue on the rate function of the LDP in
  the non-convex case. Note that the condition $\kappa_1 <2\kappa_2 -1$ is merely
  technical and yields to easier computations.
\end{rem}

\begin{proof} In order to prove Proposition \ref{second-eigen-nc}, 
we shall prove that there exists some point $z^\star$ such that 
$$
I_{\check Z}(z^\star)<\infty \qquad \textrm{while}\qquad \bar{I}_{\check Z}(z^\star)=\infty.
$$
Denote by $z=(x,y,x',y',r)$ and by ${\mathcal A}$ the convex set
$$
{\mathcal A}=\{z\in \mathbb{R}^5 ;\ x>0,\, y>0,\, x'=x,\, y'=y,\, r^2\le xy \}. 
$$
Then Cram\'er's theorem yields the LDP for $\tilde{L}_n({\check Z})$ with good rate function
$$
\Gamma^*(z) = \frac{x+y}2 -\frac 12 \log(xy -r^2) +\Delta(z\mid {\mathcal A}).
$$
Denote by ${\mathcal B}_{\kappa}$ the following non-convex set: 
$$
{\mathcal B}_{\kappa}=\{z\in \mathbb{R}^5;\ x>0,\, y>0,\, x'=\kappa x,\, y'=\kappa y,\, |r|=\sqrt{xy}\} 
$$
One can prove that $\pi_n^1({\check Z})$ and $\pi_n^2({\check Z})$ satisfy the LDP with respective rate functions
$$
I_1(z)= \frac{x+y}2 +\Delta(z\mid {\mathcal B}_{\kappa_1}) \qquad \textrm{and}
\qquad I_2(z)= \frac{x+y}2 +\Delta(z\mid {\mathcal B}_{\kappa_2}).
$$ 
The contraction principle then yields  
\begin{eqnarray*}
I_{\check Z}(z)&=&\inf_{z_0+z_1 +z_2=z}\{ \Gamma^*(z_0)+ I_1(z_1) + I_2(z_2)\}\\
\bar{I}_{\check Z}(z)&=&\inf_{z_0+z_1 =z}\{ \Gamma^*(z_0)+ I_1(z_1) \}\\
\end{eqnarray*}
Let $z^\star=(1,1,\kappa_2,\kappa_2,0)$ then we shall prove that
\begin{equation}\label{rf-different}
I_{\check Z}(z^\star)<\infty \qquad \textrm{while}\qquad \bar{I}_{\check Z}(z^\star)=\infty.
\end{equation}
This will complete the proof of Proposition \ref{second-eigen-nc}. 

In the sequel, we use the notation $z_i=(x_i,y_i,x_i',y_i',r_i)$ with
$i\in \{0,1,2\}$. From the definition of $\bar{I}_{\check Z}$, one can
easily check that $\bar{I}_{\check Z}(z^\star)$ is finite iff the
following system of equations:
\begin{equation}\label{system}
\left\{
\begin{array}{l}
x_0 + x_1 =1\\
y_0 + y_1 =1\\
x_0 +\kappa_1 x_1 = \kappa_2\\
y_0 +\kappa_1 y_1 = \kappa_2\\
x_1 y_1 < x_0 y_0
\end{array}\right.
\end{equation}
has a solution such that $x_0>0$, $y_0>0$, $x_1>0$ and $y_1>0$. From easy computations,
such a solution should satisfy
\begin{equation}\label{system-condition} 
x_0 =\frac{\kappa_1 -\kappa_2}{\kappa_1 -1}= y_0.
\end{equation}
On the other hand, the last equation of \eqref{system} implies that $(1-x_0)^2
< x_0^2$, that is $x_0>\frac 12$. As we have assumed that
$\kappa_1 <2\kappa_2 -1$, this is not compatible with \eqref{system-condition} and 
$$
\bar{I}_{\check Z}(z^\star)=\infty.
$$
We now prove that $I_{\check Z}(z^\star)<\infty$. The mere definition of $I_{\check Z}$ yields that $I_{\check Z}(z^\star)<\infty$
iff there exists a solution to the following system
\begin{equation}\label{system2}
\left\{
\begin{array}{l}
x_0 +x_1 +x_2 =1\\
y_0 +y_1 + y_2 =1\\
x_0 +\kappa_1 x_1 +\kappa_2 x_2 =\kappa_2\\
y_0 +\kappa_1 y_1 +\kappa_2 y_2 =\kappa_2\\
r_0^2 +\epsilon_1 x_1 y_1 +\epsilon_2 x_2 y_2 =0
\end{array}\right.
\end{equation}
satisfying $x_0>0$, $y_0>0$, $x_1>0$, $y_1>0$, $x_2>0$, $y_2>0$,
$\epsilon_{1,2}=\pm 1$ and $r_0^2\leq x_0y_0$.

We can easily check that this system admits the following solution:
\begin{eqnarray*}
x_0=y_0&=& \frac{\kappa_1 -\kappa_2}{\kappa_1 +\kappa_2 -2},\\ 
x_1=y_1 &=& \frac{\kappa_2 -1}{\kappa_1 +\kappa_2 -2}=x_2=y_2,\\
\epsilon_1 = - \epsilon_2 &=& -1 \qquad \textrm{and}\qquad r_0=0.
\end{eqnarray*}
Therefore, \eqref{rf-different} is proved.
\end{proof}

\subsection{Links with the spherical integral beyond the rank-one case}

When one wants to study the asymptotics of the spherical integral 
in the case  when the matrix $A_n$  
in \eqref{spherical-integral} is of finite rank larger than one,
one is led to study the Large Deviations for empirical
means which do not fulfill the convexity assumption (Assumption
(A-\ref{convex-rf})). For example, in the rank two case, the related
empirical mean to look at is given by:
$$
L_n^{(2)} = \frac 1 n \sum {\mathbf f}^{(2)}(x_i^n)\cdot Z_i,
\textrm{ with } Z_i = \left( \begin{array}{c}
X_i^2\\
Y_i^2 \\
 X_iY_i
\end{array}\right) \textrm{ and }  \mathbf{f}^{(2)}(x)=\left({\tiny \begin{array}{ccc}
1 & 0 & 0\\
0 &1 &0 \\
x &0 &0 \\
0 &x &0 \\
0 &0 &1\\
0 & 0 & x
\end{array}}\right)
$$ 
and Theorem \ref{ldp-nonconvex} applies whenever
(A-\ref{assumption-nonconvex}) is fulfilled. It is then an easy 
application of Varadhan's Lemma to get the convergence of the
spherical integrals in the rank two case (and analogously for an arbitrary
finite rank). The example studied in
Section \ref{example-non-convex-sub} supports the feeling (although in
a very indirect way) that the asymptotics of the spherical integral in
this case should depend not only on the largest eigenvalue (as 
proved in the rank-one case in \cite{GuiMai05}) but also on the second
largest eigenvalue and maybe on other ones, the number of which is
related to the rank of $A_n$.  Unfortunatelly, the very intricate
formula of the rate function associated to the LDP in the non-convex
case gives little clue on how to relate the asymptotics of the
spherical integral to the largest eigenvalues beyond the rank-one
case.

\appendix
\section{Proof of Lemma \ref{ldp-pk}}\label{proof-ldp-pk}
\begin{proof} Let $\varepsilon>0$ be fixed. Note that
  $C_{\infty}^{\mathbf{f}} \neq \emptyset$ by Assumption
  (A-\ref{Compacity}).  Since $\C^{\mathbf{f}}$ exists by
  (\ref{tight-assumption}) and is compact by (A-\ref{Compacity}),
  there exists a finite number of $m\times d$ matrices
  $(\mathbf{a}_1,\cdots, \mathbf{a}_p)$ such that
$$
\C^{\mathbf{f}}\subset \cup_{k=1}^p B(\mathbf{a}_k,\varepsilon) \quad \textrm{where}\ B(\mathbf{a}_k,\varepsilon)
=\{\mathbf{y}\in \mathbb{R}^{m\times d}, |\mathbf{y}-\mathbf{a}_k|<\varepsilon\}.
$$
>From the cover $(B(\mathbf{a}_k,\varepsilon),1\le k\le p)$, one can easily build a partition $(\Gamma_k,1\le k\le p')$ 
where $p'\le p$ with the following properties:
\begin{itemize}
\item[-] $\C^{\mathbf{f}}\subset \cup_{k=1}^{p'} \Gamma_k$,
\item[-] $\sup \{ |\mathbf{x}-\mathbf{x}'|, (\mathbf{x},\mathbf{x}')\in \Gamma_k^2\}
\le 2\varepsilon$,
\item[-] $\mathrm{int}(\Gamma_k) \cap \C^{\mathbf{f}} \neq \emptyset$
  for $1\le k\le p'$ (in particular $\mathrm{int}(\Gamma_k) \neq \emptyset$).
\end{itemize}
Let $\mathbf{b}_{k,\varepsilon}$ be an element of $\mathrm{int}(\Gamma_k) \cap \C^{\mathbf{f}}$. Denote by 
$$
\mathbf{f}^{\varepsilon}(x) =\sum_{k=1}^{p'} \mathbf{b}_{k,\varepsilon} 1_{\Gamma_k}(\mathbf{f}(x)),\quad x\in {\mathcal X}
\qquad \textrm{and}\quad {\mathcal D}^{\varepsilon} = \cap_{k=1}^p {\mathcal D}_{\mathbf{b}_{k,\varepsilon}}.
$$
We will prove in the sequel the following facts:
\begin{enumerate}
\item The partial weighted empirical mean $\pi_n^{\varepsilon}$ defined by 
$$
\pi_n^{\varepsilon}=\frac{1}{n}
\sum_{x_i^n \in C_n} \mathbf{f}^{\varepsilon}(x_i^n)\cdot Z_i
$$
satisfies the LDP with good rate function 
$\Delta^*(z\mid {\mathcal D}^{\varepsilon})=\sup\{ \la z, \lambda  \ra,\ \lambda \in {\mathcal D}^{\varepsilon}\}$.
\item The family of random variables $(\pi_n^{\varepsilon},\varepsilon>0)$ is an exponential approximation of 
$(\pi_n)$, i.e.
$$
\lim_{\varepsilon\rightarrow 0} \limsup_{n\rightarrow \infty} \frac 1n \log \mathbb{P}\{|\pi_n^{\varepsilon} -\pi_n|>\delta\}
=-\infty,\quad \forall \delta>0.
$$
\item Finally, the family $(\pi_n,n\ge 1)$ satisfies the LDP with good rate function $\Delta^*(z\mid {\mathcal D})$.
\end{enumerate}
Let us first prove fact $(1)$.
$$
\pi_n^{\varepsilon}=\frac{1}{n}
\sum_{x_i^n \in C_n} \mathbf{f}^{\varepsilon}(x_i^n)\cdot Z_i
= \frac{\mathbf{b}_{1,\,\varepsilon}}{n}\cdot \sum_{\{x_i^n, \mathbf{f}(x_i^n)\in \Gamma_1\}} Z_i +\cdots 
+ \frac{\mathbf{b}_{p',\,\varepsilon}}{n}\cdot \sum_{\{x_i^n, \mathbf{f}(x_i^n)\in \Gamma_{p'}\}} Z_i
$$
Since the sets $(\Gamma_k)$ are disjoints, the partial empirical means
$\frac 1n \sum_{x_i^n \in \mathbf{f}^{-1}(\Gamma_k)} Z_i $ are
independent. Denote by $\phi_k(n)$ the cardinality of the set
$\{x_i^n, \mathbf{f}(x_i^n)\in \Gamma_k\}$.  One has to check that
$$
\lim_{n\rightarrow \infty} \frac{\phi_k(n)}{n}=0 \quad \textrm{and}\quad \phi_k(n)\ge 1\quad \textrm{for n large enough.}
$$
Since $\phi_k(n)\le \mathrm{card}(C_n)$, the first point is proved.
Recall now that $\mathrm{int}(\Gamma_k)\cap
C^{\mathbf{f}}_{\infty}\neq 0$. Thus Condition
(\ref{tight-assumption}) yields that for $n$ large enough, there
always exist points of $C^{\mathbf{f}}_n$ that belong to
$\Gamma_k$. In particular, $\phi_k(n)\ge 1$ eventually.
Therefore, Lemma \ref{nw} yields the LDP for $\frac 1n \sum_{x_i^n \in
  \mathbf{f}^{-1}(\Gamma_k)} Z_i $ with good rate function $I(y)$.

A straightforward application of the contraction principle \cite[Theorem 4.2.1]{DemZei98} yields the LDP 
for $\pi^{\varepsilon}_n$ with good rate function
$$
\Delta^*_{\varepsilon}(z)=\inf\left\{\sum_{k=1}^{p'} \Delta^*(y_k\mid
  {\mathcal D}_Z),\ \sum_{k=1}^{p'} \mathbf{b}_{k,\varepsilon} \cdot
  y_k = z\right\}.
$$ 
We prefer the following representation  which expresses the rate function $\Delta^*_{\varepsilon}$ as an inf-convolution:
\begin{equation}
\label{convo1}
\Delta^*_{\varepsilon}(z)=\inf\left\{\sum_{k=1}^{p'} \Delta^*(z_k\mid
\mathcal D_{\mathbf{b}_{k,\varepsilon}}),\ \sum_{k=1}^{p'} z_k = z\right\}.
\end{equation}
The rate function $\Delta^*_{\varepsilon}$ is lower semi-continuous therefore  
\cite[Theorem 16.4]{Roc70} yields:
\begin{eqnarray*}
\Delta^*_{\varepsilon} & = &\sup_{\lambda \in \Rd} \left\{ \langle \lambda, z \rangle
- \sum_{1\le k \le p^{\prime}} \Delta \left(z \mid  \bar{ \mathcal
D}_{\mathbf{b}_{k,\varepsilon}}\right)\right\}\\
& = & \sup_{\lambda \in \Rd} \left\{ \langle \lambda, z \rangle
- \Delta \left(z \mid  \cap_{1\le k\le p^\prime}\bar{ \mathcal
D}_{\mathbf{b}_{k,\varepsilon}}\right)\right\}\\
& = & \Delta^* (z \mid  \cap_{1\le k\le p^\prime}\bar{\mathcal
D}_{\mathbf{b}_{k,\varepsilon}}) \stackrel{(a)}{=} \Delta^* (z \mid  \cap_{1\le k\le p^\prime} \mathcal
D_{\mathbf{b}_{k,\varepsilon}}) = \Delta^* (z \mid  \mathcal
D_{\varepsilon}).
\end{eqnarray*}
where $(a)$ follows from Proposition \ref{support-function}.
Fact (1) is proved.

Let us now prove fact (2). We have
$$
\left| \pi_n^{\varepsilon} -\pi_n\right| \le \frac 1n \sum_{x_i^n \in C_n} \left| \mathbf{f}^{\varepsilon}(x_i^n) 
-\mathbf{f}(x_i^n) \right| |Z_i|.
$$
By the definition of $\mathbf{f}^{\varepsilon}$, if $\mathbf{f}(x_i^n)\in \Gamma_k$ then 
$\mathbf{f}^{\varepsilon} (x_i^n)=\mathbf{b}_{k,\varepsilon}$
and $|\mathbf{f}(x_i^n)-\mathbf{b}_{k,\varepsilon}|\le 2\varepsilon$. Therefore 
$|\pi_n^{\varepsilon} -\pi_n|\le \frac{2\varepsilon}{n}\sum_{x_i^n \in C_n}|Z_i|$ and 
\begin{eqnarray*}
\mathbb{P}\left\{|\pi_n^{\varepsilon} -\pi_n|>\delta\right\} 
&\le & \mathbb{P}\left( \frac 1n \sum_{x_i^n \in C_n} |Z_i| > \frac{\delta}{2\varepsilon}\right)\\
&\le & \exp\left(-\frac{n\delta \kappa}{2\varepsilon} \right) \left( \esp e^{\kappa |Z_i|}\right)^{\mathrm{card}(C_n)}
\end{eqnarray*}
where $\kappa>0$ is such that $\esp e^{\kappa|Z_i|}<\infty$. Therefore 
$$
\limsup_{n\rightarrow \infty} \frac 1n \log \mathbb{P}\left\{|\pi_n^{\varepsilon} -\pi_n|>\delta\right\} 
\le -\frac{\kappa\delta}{2\varepsilon} \xrightarrow[\varepsilon\rightarrow 0]{} -\infty, 
$$
which proves the exponential equivalence. Fact (2) is proved.

We now prove fact (3). Since $(\pi_n^{\varepsilon},\varepsilon>0)$  is an exponential approximation of $\pi_n$,
Theorem 4.2.16 (a) in \cite{DemZei98} implies that $\pi_n$ satisfies a weak LDP with rate function given by:
$$
\Upsilon(z)=\sup_{\delta>0}\liminf_{\varepsilon\rightarrow 0}\inf_{z'\in B(z,\delta)} \Delta^*_{\varepsilon}(z')
\stackrel{(\star)}{=}\sup_{\delta>0}\limsup_{\varepsilon\rightarrow 0}\inf_{z'\in B(z,\delta)} \Delta^*_{\varepsilon}(z'),
$$
where $(\star)$ is a by-product of the proof of \cite[Theorem 4.2.16]{DemZei98} (see Eq. (4.2.19) for instance).
This precisely means that $\Upsilon$ is the epigraphical limit of $\Delta^*_{\varepsilon}$
(see \cite[Chapter 7]{RocWet98} for details). In order to prove that $\Upsilon=\Delta^*(\cdot \mid {\mathcal D})$,
we first note that 
\begin{equation*}
{\mathcal D}^{\varepsilon} \xrightarrow[\varepsilon\rightarrow 0]{\mathrm{pk}}\overline{\mathcal D}.
\end{equation*}
A corollary \cite[Corollary 11.35(a)]{RocWet98} of Wijsman's theorem \cite[Theorem 11.34]{RocWet98} immediatly yields:
\begin{equation}\label{epilim-cc}
  \Upsilon(z)=
  \Delta^*(z \mid \overline{\mathcal D})=\epilim_{\varepsilon\rightarrow 0} \Delta^*(z \mid {\mathcal D}^{\varepsilon}),
\end{equation}
where $\epilim$ denotes the epigraphical limit. Since $\Delta^*(z\mid
\overline{\mathcal D})= \Delta^*(z\mid{\mathcal D})$ by Proposition
\ref{support-function}, we have $\Upsilon=\Delta^*(\cdot \mid
{\mathcal D})$.  Fact (3) is thus proved and so is Lemma \ref{ldp-pk}.
\end{proof}

\bibliography{math}
\vspace{15pt}
\noindent {\sc Myl\`ene Ma\"ida},\\
Universit\'e de Paris-Sud,\\
Equipe ``Probabilit\'es-Statistiques'', b\^atiment 425\\
91405 Orsay Cedex, France.\\
e-mail: Mylene.Maida@math.u-psud.fr\\
\\
\noindent {\sc Jamal Najim},\\ 
CNRS, T\'el\'ecom Paris\\ 
46, rue Barrault, 75013 Paris, France.\\
e-mail: najim@enst.fr\\
\\
\noindent {\sc Sandrine P\'ech\'e},\\
Institut Fourier,\\ 
100 rue des Maths, BP 74\\
38402 St Martin d'Heres, France. \\
e-mail: sandrine.peche@ujf-grenoble.fr\\
\end{document}